\documentclass{article}
\usepackage{theorem,amsmath,amssymb}
\newtheorem{Th}{Theorem}[section]
\newtheorem{Le}{Lemma}[section]
\newtheorem{Co}{Corollary}[section]

\newtheorem{Rem}{Remark}[section]

\date{}
\title{On the nonstationary Stokes system in a cone: asymptotics of solutions at infinity}
\author{Vladimir Kozlov and J\"urgen Rossmann}
\begin{document}
\maketitle

\abstract{The paper deals with the Dirichlet problem for the nonstationary Stokes system in a cone.
The authors obtain existence and uniqueness results for solutions in weighted Sobolev spaces and study the asymptotics
of the solutions at infinity. \\


}

\section*{Introduction}

Although the stationary Stokes system in domains with singular boundary points is well studied (see, e.~g.,
\cite{Dauge-89,Kozlov/M/S-94,mp-83,mps-79,mr-10,sol-79}), there are only few papers dealing with the nonstationary
Stokes system in such domains. One of them is our recent paper \cite{k/r-16} 
for the problem in a 3-dimensional cone, the 2-dimensional case was studied in \cite{r-18}.
The present paper is a continuation of  \cite{k/r-16} and deals with the Dirichlet problem for the nonstationary Stokes system
\begin{eqnarray} \label{stokes1}
&&\frac{\partial u}{\partial t} - \Delta u + \nabla p = f, \quad -\nabla\cdot u = g \ \mbox{ in } K \times (0,\infty), \\ \label{stokes2}
&& u(x,t)=0 \ \mbox{ for }x\in \partial K, \ t>0, \quad u(x,0)=0 \ \mbox{ for }x\in K,
\end{eqnarray}
where $K$ is a 3-dimensional cone with vertex at the origin. In \cite{k/r-16}, we obtained existence, uniqueness and regularity results
for solutions of the problem (\ref{stokes1}), (\ref{stokes2}) in weighted Sobolev spaces, where the weight function is a power of the
distance from the vertex of the cone. However, the results in \cite{k/r-16} are not optimal, for example, in the case
that the cone $K$ is contained in a half-space. One goal (but not the main goal) of the present paper is to improve the results
for this particular case.

When considering solutions of the problem (\ref{stokes1}), (\ref{stokes2}), the question on the behavior of the solutions
both in a neighborhood of the vertex of the cone and at infinity arises. We deal here with the asymptotics
of the solution at infinity. The asymptotics near the vertex of the cone is the subject of a forthcoming paper.
Concerning the behavior at infinity, the results for the Stokes system are completely different from those for the heat equation
and other parabolic problems given in  \cite{dCN-11}, \cite{kozlov-88}--\cite{k/m-87}, \cite{k/r-11,k/r-12}. 
We show in this paper that the solution is a finite sum of singular terms and a (more regular) remainder,
where the singular terms depend on the eigenvalues of the Beltrami operator $\delta$ with Neumann boundary conditions on the intersection
$\partial \Omega$ of $\partial K$ with the unit sphere $S^2$. In the case of the heat equation and other parabolic problems,
such singular terms do not appear. It is a feature of the nonstationary Stokes system that eigenvalues
of two different operator pencils appear in solvability and regularity results. Besides the eigenvalues of the Neumann problem for the
Laplace operator, one has to consider the eigenvalues of the pencil generated by the stationary Stokes system. In analogous results for
other parabolic  problems (see \cite{kozlov-88,kozlov-89,kozlov-91}), only the eigenvalues of one operator pencil play a role.

The paper consists of three sections. Sections 1 and 2 are concerned with the parameter-depending problem
\begin{equation} \label{par1}
s\, \tilde{u} - \Delta \tilde{u} + \nabla \tilde{p} = \tilde{f}, \quad -\nabla\cdot \tilde{u} = \tilde{g}  \
   \mbox{ in } K, \quad \tilde{u} = 0 \ \mbox{ on } \partial K,
\end{equation}
which arises after the Laplace transformation with respect to the time $t$. Here, $s$ is a complex number, $\mbox{Re}\, s \ge 0$, $s\not=0$.
Section 1 deals with the solvability of this problem in weighted Sobolev spaces
In Subsections 1.3 and 1.4, we recall the main results of \cite{k/r-16}. In particular, an existence and uniqueness result
for solutions of this problem in the space $E_\beta^2(K)\times V_\beta^1(K)$ was obtained in \cite{k/r-16}.
Here $V_\beta^l(K)$ denotes the weighted Sobolev space of all functions (vector-functions) with finite norm
\begin{equation} \label{Vbeta}
\| u\|_{V_\beta^l(K)} = \Big( \int_K \sum_{|\alpha|\le l} r^{2(\beta-l+|\alpha|)}\big| \partial_x^\alpha u(x)\big|^2\, dx\Big)^{1/2},
\end{equation}
while $E_\beta^l(K)$ is the  weighted Sobolev space with the norm
\begin{equation} \label{Ebeta}
\| u\|_{E_\beta^l(K)} = \Big( \int_K \sum_{|\alpha|\le l} \big( r^{2\beta}+r^{2(\beta-l+|\alpha|)}\big)\big|
  \partial_x^\alpha u(x)\big|^2\, dx\Big)^{1/2},
\end{equation}
$r=|x|$ denotes the distance of the point $x$ from the vertex of the cone.
As was shown in \cite{k/r-16}, there are two neighboring $\beta$-intervals for which an existence and uniqueness result in the space
$E_\beta^2(K)\times V_\beta^1(K)$ holds, namely the intervals
\begin{equation} \label{interval}
\frac 12 - \lambda_1 <\beta < \frac 12 \quad\mbox{and}\quad  \frac 12 < \beta < \min\Big(\mu_2 +\frac 12\, , \lambda_1 +\frac 32\Big)
\end{equation}
Here, $\lambda_1$ and $\mu_2$ are positive numbers depending on the cone. More precisely, $\lambda_1$ is
the smallest positive eigenvalue of the operator pencil ${\cal L}(\lambda)$ generated by the Dirichlet problem for
the stationary Stokes system, while $\mu_2$ is the smallest positive eigenvalue of the operator pencil
${\cal N}(\lambda)$ generated by the Neumann problem for the Laplacian, respectively ($\mu_2(\mu_2+1)$ is the smallest positive eigenvalue of the operator
$-\delta$ with Neumamn boundary conditions, see Subsection 1.3).
The eigenvalue $\lambda_1$ is not greater than 1 since $\lambda=1$ is always an eigenvalue of the pencil ${\cal L}(\lambda)$.
In the case $\lambda_1<1$, the inequalities (\ref{interval}) for $\beta$ are sharp. However, the existence and uniqueness result
given in \cite{k/r-16} can be improved in the case that $\lambda_1$ is equal to 1 and simple. This is done in Section 1.5.
For example, the eigenvalue $\lambda=1$ is the smallest positive eigenvalue and simple if $\overline{K}\backslash \{ 0\}$ is contained
in a half-space $\alpha_1 x_1+\alpha_2 x_2+\alpha_3 x_3>0$. In this case we obtain the following weaker conditions on $\beta$,
under which an existence and uniqueness result in the space
$E_\beta^2(K)\times V_\beta^1(K)$ holds (see Theorem \ref{Bt1}):
\begin{equation} \label{interval1}
\max\big(-\mu_2-\frac 12,\, \frac 12-\mbox{Re}\, \lambda_2\big) < \beta < \min\big(\mu_2+\frac 12,\, \mbox{Re}\, \lambda_2+\frac 32\big), \ \
   \beta\not= \pm \frac 12\, , \ \ \beta\not= \frac 52\, .
\end{equation}
Here $\lambda_2$ is the eigenvalue of ${\cal L}(\lambda)$ with smallest real part $>1$. The uniqueness of the solution holds even for
$-\mu_2 -\frac 12 < \beta < \mbox{Re}\, \lambda_2 +\frac 32$ (see Lemma \ref{Bl2}).

Furthermore, we prove a regularity assertion for the
solution $(\tilde{u},\tilde{p}) \in E_\beta^2(K)\times V_\beta^1(K)$ of (\ref{par1}). If for example
\[
\tilde{f} \in E_\gamma^0(K), \ \tilde{g} \in V_\gamma^1(K)\cap(V_{-\gamma}^1(K))^*, \quad
\frac 12-\lambda_1<\beta< \frac 12 <\gamma < \min\big(\mu_2+\frac 12,\, \mbox{Re}\, \lambda_2+\frac 32\big),
\]
$\gamma\not=\frac 52$ and $\displaystyle \int_K \tilde{g}\, dx=0$, then it follows from Lemmas \ref{l1} and \ref{l2} of the present paper
that $\tilde{u} \in E_\gamma^2(K)$ and $\tilde{p}\in V_\gamma^1(K)$. This is not true, if the integral of $\tilde{g}$ over $K$ is not equal to zero.
Then we can represent $(\tilde{u},\tilde{p})$ as a sum of singular terms and a remainder $(\tilde{v},\tilde{q}) \in E_\gamma^2(K)\times V_\gamma^1(K)$
For $\gamma < \min(\frac 32\, ,\mu_2+\frac 12)$, we obtain the decomposition
\[
\big( \tilde{u},\tilde{p}\big) =  \eta\big(|s|r^2\big)\, c_1(s)\, \big( u_0^{(-1)}(x,s), p_0^{(-1)}(x)\big)+ (\tilde{v},\tilde{q})
\]
with the formulas (\ref{up0}), (\ref{coeffup0}) for $u_0^{(-1)}$, $p_0^{(-1)}$ and $c_1$,
where $\eta$ is a smooth function on $(0,\infty)$, $\eta(r)=0$ for $r<1/2$ and $\eta(r)=1$ for $r>1$. In the case $\gamma>\mu_2+\frac 12$,
additional singular terms appear, i.~e., we obtain a decomposition
\[
(\tilde{u},\tilde{p}) = \eta\big(|s|r^2\big)\ \sum_{j,k} c_{j,k}(s)\, \big( u_0^{(-j,k)}, p_0^{(-j,k)}\big) + (\tilde{v},\tilde{q}),
\]
where $\big( u_0^{(-j,k)}, p_0^{(-j,k)}\big)$ are singular functions depending on the eigenvalues of the Beltrami operator with Neumann boundary conditions
(see Theorems \ref{Ct2} and \ref{Ct3}).

In Section 3, we consider the time-dependent problem (\ref{stokes1}), (\ref{stokes2}).
The results of Section 1 enable us to obtain solvability results in weighted Sobolev spaces and regularity results for the solutions.
Partially, these results can be found in our paper \cite{k/r-16}. In the present paper, we weaken the conditions on the weight parameter $\beta$ for the
the case that $\lambda_1=1$ is the smallest positive eigenvalue of the pencil ${\cal L}(\lambda)$ and simple.
In particular, there exists a unique solution $(u,p) \in W_\beta^{2,1}(Q) \times L_2({\Bbb R}_+,V_\beta^1(K))$ for arbitrary
$f\in L_2({\Bbb R}_+,V_\beta^0(K))$, $g\in L_2({\Bbb R}_+,V_\beta^1(K))$, $\partial_t u \in L_2({\Bbb R}_+,(V_{-\beta}^1(K))^*)$
if $\frac 12-\lambda_1 < \beta <\frac 12$. Here, $W_\beta^{2,1}(Q)$ is the space of all $u\in  L_2({\Bbb R}_+,V_\beta^2(K))$ such that
$\partial_t u \in L_2({\Bbb R}_+,V_\beta^0(K))$. By means of the results of Section 3, we describe the asymptotics of this solution at infinity.
We prove that the velocity $u$ is a finite sum of terms
\[
S^{(j,k)}(x,t) = \int_0^t \int_K \big( K_u^{(j,k)}(x,y,t-\tau)\, f(y,\tau)  + H_u^{(j,k)}(x,y,t-\tau)\, g(y,\tau)\big)\, dy\, d\tau
\]
and a remainder $v\in W_\gamma^{2,1}(Q)$, $\gamma > \frac 12$, and derive point estimates for the kernels $K_u^{(j,k)}$ and $H_u^{(j,k)}$. An analogous
representation holds for the pressure $p$ (see Theorem \ref{t9}).

\section{Solvability of the parameter-depending problem}

Let $\Omega$ is a subdomain of the unit sphere $S^2$ with smooth (of class $C^{2,\alpha}$) boundary $\partial\Omega$ and let
$K = \big\{ x\in {\Bbb R}^3:\ \omega=x/|x| \in \Omega\}$ be a cone with vertex at the origin.
We consider the boundary value problem
\begin{equation} \label{par3}
s\, u - \Delta u + \nabla p = f, \ \ - \nabla\cdot u = g \ \mbox{ in }K, \quad u=0\ \mbox{ on }\partial K\backslash \{ 0\}.
\end{equation}
Here, $s$ be an arbitrary complex number, $\mbox{Re}\, s \ge 0$.
This section is concerned  with the existence and uniqueness of solutions in the space $E_\beta^2(K)\times V_\beta^1(K)$.

\subsection{Weighted Sobolev spaces on the cone}

For nonnegative integer $l$ and real $\beta$, we define the weighted Sobolev spaces $V_\beta^l(K)$ and $E_\beta^l(K)$
as the sets of all functions (or vector functions) with finite norms (\ref{Vbeta}) and (\ref{Ebeta}), respectively.
Note that the spaces $V_\beta^l(K)$ and $E_\beta^l(K)$ can be also defined as the closures of
$C_0^\infty(\overline{K}\backslash \{ 0\})$ with respect to the above norms.
Furthermore, we define $\stackrel{\circ}{V}\!{}_\beta^1(K)$ and $\stackrel{\circ}{E}\!{}_\beta^1(K)$
as the spaces of all functions $u\in V_\beta^1(K)$ and $u\in E_\beta^1(K)$, respectively, which are zero on
$\partial K\backslash \{ 0\}$. The dual spaces of $\stackrel{\circ}{V}\!{}_\beta^1(K)$ and $\stackrel{\circ}{E}\!{}_\beta^1(K)$
are denoted by $V_{-\beta}^{-1}(K)$ and $E_{-\beta}^{-1}(K)$, respectively.
Since
\[
\int_K r^{2\beta-2} \, \big| u(x)\big|^2\, dx \le c\, \int_K r^{2\beta}\, \big| \nabla u(x)\big|^2\, dx
\]
for $u\in C_0^\infty(K)$, the norm
\[
\| u \| = \Big( \int_K r^{2\beta}\, \big( |u|^2 + |\nabla u|^2\big)\, dx\Big)^{1/2}
\]
is equivalent to the $E_\beta^1(K)$-norm in $\stackrel{\circ}{E}\!{}_\beta^1(K)$.

\subsection{The operator of the problem (\ref{par3})}

Obviously, the mapping
\[
E_\beta^2(K) \times V_\beta^1(K) \ni (u,p) \to f = su-\Delta u + \nabla p \in E_\beta^0(K)
\]
is continuous for arbitrary real $\beta$ and complex $s$. Furthermore, the operator $\mbox{div}$
realizes a continuous mapping from $E_\beta^2(K) \cap \stackrel{\circ}{E}\!{}_\beta^1(K)$ into
the space
\[
X_\beta^1(K)  =  E_\beta^1(K) \cap \big(V_{-\beta}^1(K)\big)^* = V_\beta^1(K) \cap \big(V_{-\beta}^1(K)\big)^*
\]
(see \cite[Section 2.1]{k/r-16}).  We denote the operator
\begin{eqnarray} \label{1t2} \nonumber
\big( E_\beta^2(K) \cap \stackrel{\circ}{E}\!{}_{\beta}^1(K)\big) \times V_\beta^1(K) \ni (u,p)
\to  \big( su-\Delta u + \nabla p, - \nabla\cdot u\big) \in E_\beta^0(K)\times X_\beta^1(K)
\end{eqnarray}
of the problem (\ref{par3}) by $A_\beta$.
Note that that the integral of $g$ over $K$ exists if $g\in X_\beta^1(K)$, $\frac 12 < \beta <\frac 52$, and that
\[
R(A_\beta) \subset E_\beta^0(K) \times \tilde{X}_\beta^1(K) \ \mbox{ if } \frac 12 < \beta < \frac 52
\]
($R(A_\beta)$ denotes the range of the operator $A_\beta$), where
\[
\tilde{X}_\beta^1(K) = \big\{ g\in X_\beta^1(K):\ \int_K g(x)\, dx =0\big\}
\]
if $\frac 12 < \beta < \frac 52$ (see \cite[Lemma 2.12]{k/r-16}). The space $\tilde{X}_\beta^1(K)$ can be also defined
as the closure of the set all $g\in C_0^\infty(\overline{K}\backslash \{ 0\})$ satisfying the condition
\begin{equation} \label{tildeX}
\int_K g(x)\, dx =0
\end{equation}
in $X_\beta^1(K)$. However, in the cases $\beta< \frac 12$ and $\beta > \frac 52$, the following assertion is true.

\begin{Le} \label{Xl}
If $\beta<\frac 12$ or $\beta > \frac 52$, then the set of all $g\in C_0^\infty(\overline{K}\backslash \{ 0\})$ satisfying the condition
{\em (\ref{tildeX})} is dense in $X_\beta^1(K)$.
\end{Le}

P r o o f. Since the set $C_0^\infty(\overline{K}\backslash\{ 0\})$ is dense in $X_\beta^1(K)$, it suffices to show that
for every $g\in C_0^\infty(\overline{K}\backslash\{ 0\})$ there exists a sequence $(g_n) \subset C_0^\infty(\overline{K}\backslash\{ 0\})$
such that
\[
g_n \to g \mbox{ in } X_\beta^1(K)\quad\mbox{and}\quad \int_K g_n(x)\, dx =0 \ \mbox{ for all }n
\]
if $\beta> \frac 52$ or $\beta< \frac 12$.
Let $\zeta$ be a differentiable function on $\overline{K}$ sucht that $\zeta(x)=0$ for $|x|<1$ and $|x|>2$ and
$\displaystyle \int_K \zeta(x)\, dx =1$. Furthermore, we set $\zeta_n(x) = n^3\, \zeta(nx)$. Then
$\int_K \zeta_n(x)\, dx=1$ and
\[
\| \zeta_n \|_{X_\beta^1(K)} \le \| \zeta_n\|_{V_\beta^1(K)} + \| \zeta_n\|_{V_{\beta+1}^0(K)} = n^{-\beta+5/2}\, \| \zeta\|_{V_\beta^1(K)}
  + n^{-\beta+1/2} \| \zeta_n\|_{V_{\beta+1}^0(K)}.
\]
This means that $\zeta_n\to 0$ in $X_\beta^1(K)$ as $n \to \infty$ if $\beta> \frac 52$.
Consequently, the sequence of the functions
\[
g_n(x) = g(x) - \zeta_n(x)\, \int_K g(x)\, dx
\]
converges to $g$ in $X_\beta^1(K)$ if $\beta>\frac 52$. Analogously, it can be shown that the sequence of the functions
\[
h_n(x) = g(x) - n^{-3}\, \zeta(n^{-1}x)\, \int_K g(x)\, dx
\]
converges to $g$ in $X_\beta^1(K)$ if $\beta<\frac 12$. Furthermore, the integrals of $g_n$ and $h_n$ over $K$ are zero. This
proves the lemma. \hfill $\Box$ \\

Let $A^*_\beta$ denote the adjoint operator of $A_\beta$. This operator is defined as a continuous mapping
\[
E_{-\beta}^0(K) \times \big( X_\beta^1(K)\big)^* \to \big( E_\beta^2(K)\big)^* \times \big( V_\beta^1(K)\big)^*.
\]
Here, $\big( X_\beta^1(K)\big)^* = V_{-\beta}^1(K) +\big( V_\beta^1(K)\big)^*$. Since the constant function
$g=1$ is an element of the space $\big( X_\beta^1(K)\big)^*$ for $\frac 12 < \beta < \frac 52$, the kernel
of $A_\beta^*$ contains the pair $(v,q)=(0,c)$ with constant $c$ in the case $\frac 12 < \beta < \frac 52$.
Other constant elements of $\mbox{ker}\, A_\beta^*$ do not exist. The kernel of the operator $A_\beta$
contains no constant elements except $(u,p)=(0,0)$.

\subsection{Normal solvability of the operator $A_\beta$}

We introduce the following operator pencils ${\cal L}(\lambda)$ and ${\cal N}(\lambda)$ generated by
the Dirichlet problem for the stationary Stokes system and the Neumann problem for the Laplacian
in the cone $K$, respectively. For every complex $\lambda$, we define the operator ${\cal L}(\lambda)$ as the mapping
\begin{eqnarray*}
&& \stackrel{\circ}{W}\!{}^1(\Omega)\times L_2(\Omega) \ni \left( \begin{array}{c} U \\ P \end{array}\right) \\ && \qquad
  \to \left( \begin{array}{c} r^{2-\lambda}\big(-\Delta r^{\lambda}U(\omega)+\nabla r^{\lambda-1}P(\omega)\big)\\*[1ex]
  -r^{1-\lambda}\nabla\cdot\big( r^\lambda U(\omega)\big)\end{array}\right) \in W^{-1}(\Omega)\times L_2(\Omega),
\end{eqnarray*}
where $r=|x|$ and $\omega=x/|x|$. The properties of the pencil ${\cal L}$ are studied, e.g., in \cite{kmr-01}.
In particular, it is known that the numbers $\lambda$, $\bar{\lambda}$ and $-1-\lambda$ are simultaneously eigenvalues
of the pencil ${\cal L}(\lambda)$ or not. The eigenvalues in the strip $-2\le \mbox{Re}\, \lambda \le 1$ are real,
and the numbers $1$ and $-2$ are always eigenvalues of the pencil ${\cal L}(\lambda)$. If $\Omega$ is contained in a half-sphere,
then $\lambda=1$ and $\lambda=-2$ are the only eigenvalues in the interval $[-2,1]$ (cf. \cite[Theorem 5.5.5]{kmr-01}).
We denote the eigenvalues with positive real part by $\lambda_j$, $j=1,2,\ldots$, while $\lambda_{-j}=-1-\lambda_j$
are the eigenvalues with negative real part,
\[
\cdots \le \mbox{Re}\, \lambda_{-2} < \lambda_{-1} < -1 < 0 < \lambda_1 < \mbox{Re}\, \lambda_2 \le \cdots
\]
Here, $0<\lambda_1\le 1$ and $-2\le \lambda_{-1}<-1$. Note that the eigenvalues $\lambda_j$ and $\lambda_{-j}$ have the same
geometric and algebraic multiplicities.

The operator ${\cal N}(\lambda)$ is defined as
\[
{\cal N}(\lambda)\, U = \Big( -\delta U-\lambda(\lambda+1)U\, , \, \frac{\partial U}{\partial\vec{n}}\big|_{\partial\Omega}\Big)\
  \quad\mbox{for } U\in W^2(\Omega).
\]
As is known (see e.g. \cite[Section 2.3]{kmr-01}), the eigenvalues of this pencil are real, and
generalized eigenfunctions do not exist. The spectrum contains, in particular, the simple eigenvalues $\mu_1=0$ and $\mu_{-1}=-1$ with the
eigenfunction $\phi_1 = const.$ The interval $(-1,0)$ is free of eigenvalues. Let $\mu_j$, $j=1,2,\ldots$, be the nonnegative eigenvalues,
and let $\mu_{-j}=-1-\mu_j$ be the negative eigenvalues of the pencil ${\cal N}(\lambda)$,
\[
\cdots < \mu_{-2} < -1 =\mu_{-1} < \mu_1 =0 < \mu_2 < \cdots.
\]
Obviously, $\mu_j$ and $\mu_{-j}$ are the solutions of the equation $\lambda(\lambda+1)=-M_j$, where $M_j$ is the $j$th eigenvalue
of the operator $-\delta$ with Neumann boundary condition on $\partial\Omega$. For the following theorem, we refer to \cite[Theorem 2.1]{k/r-16}

\begin{Th} \label{t1}
Suppose that $\mbox{\em Re}\, s \ge 0$, $|s|=1$, that the line $\mbox{\em Re}\, \lambda = -\beta+1/2$ does not contain
eigenvalues of the pencil ${\cal L}(\lambda)$, and that $-\beta-1/2$ is not an eigenvalue of the pencil ${\cal N}(\lambda)$.
Then the range of the operator {\em (\ref{1t2})} is closed and the kernel has finite dimension.
\end{Th}

Note that the condition on the eigenvalues of the pencils ${\cal L}(\lambda)$ and ${\cal N}(\lambda)$  in Theorem \ref{t1}
is necessary (cf. \cite[Lemmas 2.7, 2.8]{k/r-16}). In \cite{k/r-16} it was also shown that the following regularity assertion
for solutions of the problem (\ref{par3}) is true.

\begin{Le}  \label{l5b}
Suppose that $(u,p) \in E_\beta^2(K) \times V_\beta^1(K)$ is a solution of the problem {\em (\ref{par3})}, where $\mbox{\em Re}\, s\ge 0$
and $|s|=1$, $f\in E_\beta^0(K) \cap E_\gamma^0(K)$ and $g\in X_\beta^1(K)\cap X_\gamma^1(K)$.
We assume that one of the following two conditions is satisfied:
\begin{itemize}
\item[{\em (i)}] $\beta<\gamma$ and the interval $-\gamma-1/2 \le \lambda \le -\beta-1/2$ does not contain eigenvalues of the pencil
${\cal N}(\lambda)$,
\item[{\em (ii)}] $\beta>\gamma$ and the strip $-\beta+1/2 \le \mbox{\em Re}\, \lambda\le -\gamma+1/2$ is free of eigenvalues
  of the pencil ${\cal L}(\lambda)$.
\end{itemize}
Then $u\in E_\gamma^2(K)$, $p\in V_\gamma^1(K)$ and
\[
\| u\|_{E_\gamma^2(K)} +\| p\|_{V_\gamma^1(K)} \le c\, \Big( \| f\|_{E_\gamma^0(K)} + \| g\|_{X_\gamma^1(K)}
  + \| u\|_{E_\beta^2(K)} +\| p\|_{V_\beta^1(K)}\Big).
\]
Here, the constant $c$ is independent of $f$, $g$ and $s$.
\end{Le}

\subsection{Bijectivity of the operator $A_\beta$}

The following lemma is essentially proved in \cite{k/r-16}.

\begin{Le} \label{l6}
Suppose that $\mbox{\em Re}\, s \ge 0$, $s\not=0$ and  $-\mu_2 -1/2 < \beta < \lambda_1 +3/2$.
Then $A_\beta$ is injective.
\end{Le}

P r o o f.
By \cite[Lemma 2.10]{k/r-16}, the operator $A_\beta$ is injective if
$-\mu_2 -1/2 < \beta < \lambda_1 +3/2$ and $\beta\not=-1/2$. We prove the injectivity for $\beta=-1/2$.
Let $\zeta$ be a smooth function with compact support which is equal to one near the vertex of the cone $K$, and let
$\eta=1-\zeta$. Furthermore, let $\varepsilon$ be a sufficiently small positive number.
Suppose that $(u,p) \in \mbox{ker}\, A_{-1/2}$. Then $\eta p \in V_{-\varepsilon-1/2}^1(K)$  and
\[
\int_K \nabla(\eta p)\cdot \nabla q\, dx = \langle F,q\rangle,
\]
where
\[
\langle F,q\rangle = \int_K \Delta u\cdot \nabla(\eta q)\, dx + \int_K \nabla\eta \cdot (p\nabla q - q\nabla p)\, dx
\]
for all $q\in V_{\varepsilon+1/2}^1(K)$. By \cite[Lemma 2.5]{k/r-16}, the functional $F$ is continuous both on
$V_{\varepsilon+1/2}^1(K)$ and on $V_{\varepsilon+1}^2(K)$ if $0<\varepsilon\le \frac 12$.
Since the interval $\varepsilon -\frac 12 \le \lambda \le \varepsilon$ contains only the simple eigenvalue
$\lambda=0$ of the pencil ${\cal N}(\lambda)$ for small positive $\varepsilon$, it follows from \cite[Lemma 2.6]{k/r-16} that
$\eta p = c + q_0$, where $q_0\in V_{-1-\varepsilon}^0(K)$ and $c$ is a constant. Consequently, $p= c+q_1$,
where $q_1=q_0+\zeta p \in V_{-1-\varepsilon}^1(K)$. Furthermore, $u\in E_{-1/2}^2(K) \subset V_{-1-\varepsilon}^1(K)$ for
$\varepsilon \le \frac 12$. Since the pair $(u,q_1)$ is also a solution of the Dirichlet problem for the system
\[
(s-\Delta)\, u + \nabla q_1 = 0, \quad \nabla \cdot u =0\ \mbox{ in }K,
\]
it follows from \cite[Lemma 2.4]{k/r-16} that $u\in E_{-\varepsilon}^1(K)$, $q_1\in V_{-\varepsilon}^1(K)$ and, consequently,
$(u,q_1) \in \mbox{ker}\, A_{-\varepsilon}$. This means that $u=0$ and $q_1=0$, i.~.e., $p=c$. However, the space
$V_{-1/2}^1(K)$ contains no constants except $c=0$. Hence, the kernel of $A_{-1/2}$ is trivial. The proof is complete.
\hfill $\Box$ \\

Furthermore, the following result was proved in \cite[Theorems 2.3--2.5]{k/r-16}.

\begin{Th} \label{t3}
Suppose that $\mbox{\em Re}\, s \ge 0$ and $s\not=0$.

{\em 1)} If $\frac 12 -\lambda_1 <\beta < \frac 12$, then the operator $A_\beta$ is an isomorphism onto $E_\beta^0(K) \times X_\beta^1(K)$,
and the estimate
\begin{equation} \label{1t3}
\| u\|_{V_\beta^2(K)} + |s|\, \| u\|_{V_\beta^0(K)} + \| p\|_{V_\beta^1(K)} \le c\, \Big( \| f\|_{V_\beta^0(K)}
  + \| g \|_{V_\beta^1(K)} + |s|\, \| g\|_{(V_{-\beta}^1(K))^*}\Big)
\end{equation}
is valid for every solution $(u,p)\in E_\beta^2(K)\times V_\beta^1(K)$ of the problem {\em (\ref{par3})}.

{\em 2)} If $\frac 12 <\beta < \min\big( \mu_2 +\frac 12\, , \, \lambda_1 +\frac 32\big)$, then
the operator $(u,p) \to (f,g)$ of the problem {\em (\ref{par3})} is an isomorphism onto $E_\beta^0(K) \times \tilde{X}_\beta^1(K)$,
and the estimate {\em (\ref{1t3})} is valid for every solution
$(u,p)\in E_\beta^2(K)\times V_\beta^1(K)$ of the problem {\em (\ref{par3})}.
\end{Th}

The condition $\frac 12 -\lambda_1 <\beta < \min\big( \mu_2 +\frac 12\, , \, \lambda_1 +\frac 32\big)$, $\beta\not=\frac 12$
for the bijectivity of the operator $A_\beta$ is sharp if $\lambda_1 <1$ (cf. \cite[Lemmas 2.14, 2.15, 2.17]{k/r-16}).
Furthermore, the following regularity assertion can be easily deduced from Theorem \ref{t3} and \cite[Lemma 2.13]{k/r-16}.

\begin{Le} \label{l1}
Let $(u,p) \in E_\beta^2(K)\times V_\beta^1(K)$ be a solution of the problem {\em (\ref{par3})}, where
\[
f\in E_\beta^0(K)\cap E_\gamma^0(K), \quad g\in X_\beta^1(K)\cap X_\gamma^1(K),
\]
$\frac 12 -\lambda_1 < \beta,\gamma <\frac 12 +\min (\mu_2, \lambda_1+1)$, $\beta,\gamma \not= \frac 12$, $s\not=0$ and $\mbox{\em Re}\, s\ge 0$.
In the case $\max(\beta,\gamma)>\frac 12$, we assume in addition that $g$ satisfies the condition {\em (\ref{tildeX})}. Then
$u\in E_\gamma^2(K)$, $p\in V_\gamma^1(K)$.
\end{Le}

\subsection{The case that the eigenvalue $\lambda_1$ is equal to 1 and simple}

As was mentioned above, the number $\lambda=1$ is always an eigenvalue of the pencil ${\cal L}(\lambda)$ with the corresponding
constant eigenvector $(0,1)$. In the case that $\overline{\Omega}$ is contained in a half-sphere, the eigenvalues $\lambda=1$
and $\lambda=-2$ are simple (have geometric and algebraic multiplicity 1), and all other eigenvalues lie outside the strip
$-2 \le \mbox{Re}\, \lambda \le 1$.

We assume in this subsection that (as in the just described case), the smallest positive eigenvalue of the pencil
${\cal L}(\lambda)$ is $\lambda_1=1$ and that this eigenvalue is simple. Since the strip $-2\le \mbox{Re}\, \lambda
\le 1$ contains only real eigenvalues of the pencil ${\cal L}(\lambda)$ (cf. \cite[Theorem 5.3.1]{kmr-01}), it follows then that
\[
\cdots\le  \mbox{Re}\, \lambda_{-2} < -2 =\lambda_{-1} < \lambda_1 = 1 < \mbox{Re}\, \lambda_2 \le \cdots .
\]
Then the results of Lemma \ref{l6} and Theorem \ref{t3} can be improved.

\begin{Le} \label{Bl1}
Suppose that $\mbox{\em Re}\, s\ge 0$, $s\not=0$, $\lambda_1=1$ and that $\lambda_1$ is a simple eigenvalue of the pencil ${\cal L}(\lambda)$.
Then $A_\beta$ is an isomorphism onto $E_\beta^0(K) \times X_\beta^1(K)$ if $\max(-\mu_2 - \frac 12\, ,\, \frac 12-\mbox{\em Re}\, \lambda_2)
<\beta<-\frac 12$.
\end{Le}

P r o o f.
Suppose that $\max(-\mu_2 - \frac 12\, ,\, \frac 12-\mbox{Re}\, \lambda_2)<\beta<-\frac 12$. Then the kernel of $A_\beta$
is trivial (see Lemma \ref{l6}) and the range of $A_\beta$ is closed (see Theorem \ref{t1}). Thus, it suffices to show
that the problem (\ref{par3}) is solvable in $E_\beta^2(K) \times V_\beta^1(K)$ for arbitrary
$f\in C_0^\infty(\overline{K}\backslash \{ 0\})$, $g\in C_0^\infty(\overline{K}\backslash \{ 0\})$.
By Theorem \ref{t3}, there exists a solution $(u,p)\in E_\gamma^2(K)\times V_\gamma^1(K)$, where $-\frac 12 < \gamma <\frac 12$.
We show that $u\in E_\beta^2(K)$ and $p-c\in V_\beta^1(K)$ for a certain constant $c$.
Let $\zeta$ be a two times continuously differentiable function with compact support in $\overline{K}$ which is equal to one
in a neighborhood of the vertex of the cone $K$, and let $\eta=1-\zeta$. Then $\zeta(u,p) \in V_\gamma^2(K) \times
V_\gamma^1(K)$, $-\nabla \cdot (\zeta u) = \zeta g - u\cdot \nabla\zeta \in V_\beta^1(K)$ and
\[
-\Delta(\zeta u)+\nabla(\zeta p) = F, \ \mbox{where }F=\zeta f -s\zeta u + \zeta \Delta u - \Delta(\zeta u) +p\nabla\zeta \in V_\gamma^0(K) \cap V_{\gamma-2}^0(K)
\]
If $\gamma-2\le \beta\le \gamma$, then $V_\gamma^0(K) \cap V_{\gamma-2}^0(K) \subset V_\beta^0(K)$, and we conclude
from well-known regularity results for solutions of elliptic boundary value problems (see, e.~g., \cite[Chapter 3, Theorem 5.5]{np-94})
that $\zeta u\in V_\beta^2(K)$ and $\zeta p-c\in V_\beta^0(K)$ since $\lambda_1=1$ is the only eigenvalue of the pencil ${\cal L}(\lambda)$
in the strip $\frac 12-\gamma \le \mbox{Re}\, \lambda \le \frac 12 - \beta$.
Obviously, $\eta(u,p) \in E_\beta^2(K)\times V_\beta^1(K)$. Thus, $u\in E_\beta^2(K)$ and $p-c \in V_\beta^1(K)$.
If $\beta < \gamma-2$, then we conclude first that $\zeta u \in V_{\gamma-2}^2(K)$ and $\zeta p -c \in V_{\gamma-2}^0(K)$.
In this case, we conclude that $F \in V_{\gamma-4}^0(K) \cap V_{\gamma-2}^0(K)$.
Applying again the regularity result in \cite[Chapter 3, Theorem 5.5]{np-94}, we obtain $\zeta u\in V_\beta^2(K)$, $\zeta p -c
\in V_\beta^1(K)$ if $\gamma-4\le \beta \le \gamma$. In this way, after finitely many steps, we get
$u\in E_\beta^2(K)$ and $p-c\in V_\beta^1(K)$ if $\max(-\mu_2 - \frac 12\, ,\, \frac 12-\mbox{Re}\, \lambda_2)<\beta<-\frac 12$.
Obviously, the pair $(u,p-c)$ is also a solution of the problem (\ref{par3}). Thus, it is shown that (\ref{par3}) is solvable in
the space $E_\beta^2(K) \times V_\beta^1(K)$ for arbitrary $f\in C_0^\infty(\overline{K}\backslash \{ 0\})$
and $g\in C_0^\infty(\overline{K}\backslash \{ 0\})$. The proof of the lemma is complete. \hfill $\Box$ \\

The last lemma allows us to improve the result of Lemma \ref{l6} if $\lambda_1 =1$ and $\lambda_1$ is a simple eigenvalue.

\begin{Le} \label{Bl2}
Suppose that $\mbox{\em Re}\, s\ge 0$, $s\not=0$, $\lambda_1=1$ and that $\lambda_1$ is a simple eigenvalue of the pencil ${\cal L}(\lambda)$.
Furthermore, we assume that $-\mu_2 -\frac 12 < \beta < \mbox{\em Re}\, \lambda_2 +\frac 32$. Then the operator
$A_\beta$  is injective.
\end{Le}

P r o o f.
By Lemma \ref{l6}, the operator $A_\beta$ is injective for $-\mu_2 -\frac 12 < \beta < \frac 52$.
We show that $A_\beta$ is injective for $\frac 52 \le \beta < \delta +\frac 52$, where $\delta= \min(\mu_2,\mbox{Re}\, \lambda_2-1)$.
Suppose that the kernel of $A_\beta$ is not trivial for one $\beta$ in the interval $\frac 52 \le \beta < \delta+\frac 52$. As was shown
in \cite{k/r-16} (see Formula (34) in \cite{k/r-16}), there is the relation
\[
\mbox{ker}\, A_\gamma^* \supset \mbox{ker}\, A_\beta \ \mbox{ if } -\beta \le \gamma \le 2-\beta.
\]
Thus, the kernel of $A^*_\gamma$ is not trivial for $-\beta\le \gamma \le 2-\beta$. However the interval
$[-\beta,2-\beta]$ has a nonempty intersection with the interval $(\max(-\mu_2 - \frac 12\, ,\, \frac 12-\mbox{Re}\, \lambda_2),-\frac 12)$
since $2-\beta>-\mu_2 - \frac 12$ and $2-\beta>\frac 12-\mbox{Re}\, \lambda_2$ for $\beta<\delta+\frac 52$. This contradicts Lemma \ref{Bl1}.
Consequently, the kernel of $A_\beta$ is trivial for  $\frac 52 \le \beta < \delta+\frac 52$. Furthermore, it follows from Lemma \ref{l5b} that
$\mbox{ker}\, A_\beta \subset \mbox{ker}\, A_\gamma=\{ 0\}$ if $\frac 52 <  \gamma < \delta + \frac 52\le \beta< \mbox{Re}\, \lambda_2 +\frac 32$,
i.~e., $\mbox{Re}\, \lambda_{-2} < \frac 12 - \beta < \frac 12 -\gamma < -2 = \lambda_{-1}$.
This proves the lemma. \hfill $\Box$ \\

In the case $\mu_2>2$, the following lemma holds.

\begin{Le} \label{Bl3}
Suppose that $\mbox{\em Re}\, s\ge 0$, $s\not=0$, $\lambda_1=1$ and that $\lambda_1$ is a simple eigenvalue of the pencil ${\cal L}(\lambda)$.
Furthermore, we assume that $\mu_2 >2$ and $\frac 52 < \beta < \min(\mu_2 + \frac 12,\mbox{\em Re}\, \lambda_2 + \frac 32)$. Then
$A_\beta$ is an isomorphism onto $E_\beta^0(K) \times X_\beta^1(K)$.
\end{Le}

P r o o f.
By Lemma \ref{Bl2}, the operator $A_\beta$ is injective. Furthermore, under the given condition on $\beta$, the  number $-\frac 12-\beta$
is not an eigenvalue of the pencil ${\cal N}(\lambda)$ (since $\mu_{-2}=-\mu_2-1 < -\frac 12 -\beta < -3 < \mu_{-1}=-1$), and the line
$\mbox{Re}\, \lambda = \frac 12 -\beta$ is free of eigenvalues of the pencil ${\cal L}(\lambda)$ (since $\mbox{Re}\, \lambda_{-2} =
-\mbox{Re}\lambda_2 -1< \frac 12 -\beta <-2 =\lambda_{-1}$). Hence, the range of $A_\beta$ is closed (see Theorem \ref{t1}).
Let $f\in C_0^\infty(\overline{K}\backslash \{ 0\})$, $g\in C_0^\infty(\overline{K}\backslash \{ 0\})$, and let $g$ satisfy the condition
(\ref{tildeX}). Then by Theorem \ref{t3}, there exists a solution $(u,p)\in E_\gamma^2(K) \times V_\gamma^1(K)$, where $\frac 12 < \gamma < \frac 52$.
Since $\mu_{-2} =-1-\mu_2 < -\frac 12 - \beta < -\frac 12 - \gamma < -1=\mu_{-1}$, the interval $-\frac 12 -\beta \le \lambda \le -\frac 12 -\gamma$
does not contain eigenvalues of the pencil ${\cal N}(\lambda)$. Therefore, it follows from Lemma \ref{l5b} that $u\in E_\beta^2(K)$ and
$p\in V_\beta^1(K)$. Consequently, the range of the operator $A_\beta$ contains the set of all pairs $(f,g)$
of $C_0^\infty(\overline{K}\backslash \{ 0\})$-functions satisfying the condition (\ref{tildeX}).
By Lemma \ref{Xl}, the set of all $g\in C_0^\infty(\overline{K}\backslash \{ 0\})$ satisfying the condition (\ref{tildeX}) is dense
in $X_\beta^1(K)$. Thus, the range of the operator $A_\beta$ is the set $E_\beta^0(K)\times X_\beta^1(K)$. The proof is complete. \hfill $\Box$. \\

We can summarize the results of Theorem \ref{t3} and Lemmas \ref{Bl1}, \ref{Bl3} as follows.

\begin{Th} \label{Bt1}
Suppose that $\mbox{\em Re}\, s\ge 0$, $s\not=0$, $\lambda_1=1$ and that $\lambda_1$ is a simple eigenvalue of the pencil ${\cal L}(\lambda)$.
Then the following assertions are true.

{\em 1)} If $\max(-\mu_2-\frac 12\, ,\, \frac 12 -\mbox{\em Re}\, \lambda_2)< \beta <\frac 12$ and $\beta\not=-\frac 12$, then
  $A_\beta$ is an isomorphism onto the space $E_\beta^0(K) \times X_\beta^1(K)$.

{\em 2)} If $\frac 12 < \beta < \min(\mu_2 +\frac 12\, ,\, \frac 52)$, then
  $A_\beta$ is an isomorphism onto $E_\beta^0(K) \times \tilde{X}_\beta^1(K)$.

{\em 3)} If $\mu_2 > 2$ and $\frac 52 < \beta < \min(\mu_2 +\frac 12\, ,\, \mbox{\em Re}\, \lambda_2 +\frac 32)$, then
  $A_\beta$ is an isomorphism onto $E_\beta^0(K) \times X_\beta^1(K)$.
\end{Th}

Note that the operator $A_\beta$ is not Fredholm for the values $\pm \frac 12$ and $\frac 52$ of $\beta$.
Indeed, for $\beta=-\frac 12$, we have $\frac 12 - \beta=\lambda_1$ and $-\frac 12 -\beta = \mu_1$. If $\beta=\frac 12$, then
$-\frac 12 -\beta=\mu_{-1}$, while $\frac 12 - \beta = \lambda_{-1}$ if $\beta=\frac 52$. Thus, the conditions of
Theorem \ref{t1} on $\beta$ are not satisfied for $\beta=\pm \frac 12$ and $\beta=\frac 52$. The following lemma
shows that the bounds for $\beta$ in Theorem \ref{Bt1} are sharp.

\begin{Le} \label{Bl4}
Suppose that $\mbox{\em Re}\, s\ge 0$, $s\not=0$, $\lambda_1=1$ and that $\lambda_1$ is a simple eigenvalue of the pencil ${\cal L}(\lambda)$.
Then the following assertions are true.

{\em 1)} If $\mu_2>\mbox{\em Re}\, \lambda_2-1$ and $\max(-\mu_2-\frac 12\, ,\,  -\mbox{\em Re}\, \lambda_2 -\frac 32)< \beta <-\mbox{\em Re}\, \lambda_2 +\frac 12$,
then the kernel of  $A_\beta^*$ is not trivial.

{\em 2)} If $\mu_2<\mbox{\em Re}\, \lambda_2-1$ and $-\mbox{\em Re}\, \lambda_2 +\frac 12< \beta <-\mu_2 -\frac 12$,
then the kernel of  $A_\beta$ is not trivial.

{\em 3)} If $\mu_2<\mbox{\em Re}\, \lambda_2+1$ and $\mu_2+\frac 12< \beta <\min( \mu_2+\frac 52 \, ,\, \mbox{\em Re}\, \lambda_2 +\frac 32)$,
then the kernel of  $A_\beta^*$ contains nonconstant elements.

{\em 4)} If $\mu_2>\mbox{\em Re}\, \lambda_2+1$ and $\mbox{\em Re}\, \lambda_2 +\frac 32 < \beta < \mu_2+\frac 12$,
then the kernel of  $A_\beta$ is not trivial.
\end{Le}

P r o o f. 1) Let $\zeta,\eta$ be the same smooth functions as in the proof of Lemma \ref{Bl1}. Furthermore, let $u=\zeta\, r^{\lambda_2}\phi(\omega)$,
$p=\zeta r^{\lambda_2-1}\psi(\omega)$, where $(\phi,\psi)$ is an eigenvector of the pencil ${\cal L}$ corresponding to $\lambda_2$. Then
$(u,p) \in E_\gamma^2(K)\times V_\gamma^1(K)$ for arbitrary $\gamma > -\mbox{Re}\, \lambda_2 +\frac 12$, but $(u,p) \not\in
E_\beta^2(K)\times V_\beta^1(K)$. Moreover, $(s-\Delta)u + \nabla p \in E_\beta^2(K)$ and $\nabla\cdot u \in X_\beta^1(K)$ since $\beta>
-\mbox{Re}\, \lambda_2 -\frac 32$. We assume that the kernel of $A_\beta^*$ is trivial and, consequently, $A_\beta$ is an isomorphism onto
$E_\beta^0(K)\times X_\beta^1(K)$ (see Lemma \ref{Bl2}). Then there exists a pair $(v,q)\in E_\beta^2(K)\times V_\beta^1(K)$ such that
\begin{equation} \label{2Bl4}
(s-\Delta)(u-v)+\nabla(p-q)=0, \ \ \nabla\cdot(u-v)=0 \, \mbox{ in K}, \quad u=v=0\ \mbox{ on }\partial K\backslash\{ 0\}.
\end{equation}
Suppose that $-\mbox{Re}\, \lambda_2 +\frac 12 < \gamma < -\frac 12$. Then $\mu_1=0 < -\gamma-\frac 12 < -\beta-\frac 12 < \mu_2$.
Thus, we conclude from Lemma \ref{l5b} that $(v,q) \in E_\gamma^2(K)\times V_\gamma^1(K)$. This means that
$(u-v,p-q) \in \mbox{ker}\, A_\gamma$. Since $\mbox{ker}\, A_\gamma=\{ 0\}$, it follows that $(u,p)=(v,q)\in E_\beta^2(K)\times V_\beta^1(K)$.
Hence, our assumption $\mbox{ker}\, A_\beta^*=\{ 0\}$ led to a contradiction.

2) Let $\phi$ be an eigenfunction of the pencil ${\cal N}$ corresponding to the eigenvalue $\mu_2$. In the proof of
\cite[Lemma 2.17]{k/r-16} we constructed a vector function $(U,P)$ with the leading terms $p_0=r^{\mu_2}\phi(\omega)$ and $u_0=-s^{-1}\nabla p_0$
which has the following properties:
\begin{eqnarray*}
&& \eta (U,P) \in E_\beta^2(K)\times V_\beta^1(K) \ \mbox{ for }\beta<-\mu_2-\frac 12, \quad U=0\ \mbox{ on }\partial K\backslash\{ 0\}, \\
&& (s-\Delta)\, (\eta U) + \nabla(\eta P) \in E_\gamma^0(K), \ \ \nabla\cdot(\eta U)\in X_\gamma^1(K)\ \mbox{ for }\gamma<-\mu_2+\frac 32 \, .
\end{eqnarray*}
(By Lemma \ref{Cl4}, the vector functions $(U_N,P_N)$ constructed in Lemma \ref{Cl3} have these properties for $\mu=\mu_2$ and $N\ge 1$.)
However, $\eta P\not\in V_\gamma^1(K)$ for $\gamma >-\mu_2-\frac 12$. Let $-\mu_2-\frac 12 < \gamma< \min(-\mu_2+\frac 32\, ,\, -\frac 12)$. By
Theorem \ref{Bt1}, there exists a vector function $(v,q)\in E_\gamma^2(K)\times V_\gamma^1(K)$ such that
\begin{equation} \label{1Bl4}
(s-\Delta)(v-\eta U)+\nabla(q-\eta P)=0, \ \ \nabla\cdot(v-\eta U)=0 \, \mbox{ in K}, \quad v=0\ \mbox{ on }\partial K\backslash\{ 0\}.
\end{equation}
Since $\lambda_1=1 < -\gamma+\frac 12 < -\beta+\frac 12 <\mbox{Re}\, \lambda_2$, we conclude from Lemma \ref{l5b} that
$(v,q)\in E_\beta^2(K)\times V_\beta^1(K)$. Hence, $(v-\eta U,q-\eta P)$ is a nonzero element of $\mbox{ker}\, A_\beta$.

3) Let $(U,P)$ be the same vector function as in the second part of the proof, and let
$\max( -\mu_2-\frac 12\, , -\mbox{Re}\, \lambda_2 +\frac 12) < \gamma < \min(\frac 12\, ,-\mu_2+\frac 32)$, $\gamma\not=-\frac 12$. Then
\[
\eta (U,P) \in E_{-\beta}^2(K)\times V_{-\beta}^1(K) \subset E_{-\beta}^0(K) \times \big( X_\beta^1(K)\big)^*, \quad U=0\ \mbox{ on }\partial K\backslash\{ 0\}.
\]
and
\[
(s-\Delta)\, (\eta U) + \nabla(\eta P) \in E_\gamma^0(K), \ \ \nabla\cdot(\eta U)\in X_\gamma^1(K).
\]
By Theorem \ref{Bt1}, there exists a vector function $(v,q)\in E_\gamma^2(K)\times V_\gamma^1(K)$ satisfying (\ref{1Bl4}).
For $-\beta \le \gamma \le 2-\beta$, the imbeddings
\[
E_\gamma^2(K) \subset E_{-\beta}^0(K), \quad V_\gamma^1(K) \subset \big( X_\beta^1(K)\big)^*
\]
hold. Since the intervals $[-\beta,2-\beta]$ and $\max( -\mu_2-\frac 12\, -\mbox{Re}\, \lambda_2 +\frac 12) < \gamma < \min(\frac 12\, ,-\mu_2+\frac 32)$
have a nonempty intersection for $\mu_2+\frac 12< \beta <\min( \mu_2+\frac 52 \, ,\, \mbox{Re}\, \lambda_2 +\frac 32)$,
we can choose $\gamma$ such that $(v,q) \in E_{-\beta}^0(K) \times \big( X_\beta^1(K)\big)^*$. Then
$(v-\eta U,q-\eta P)$ is a nonconstant element of $\mbox{ker}\, A_\beta^*$.

4) Let $u=\zeta\, r^{-1-\lambda_2}\phi(\omega)$,
$p=\zeta r^{-2-\lambda_2}\psi(\omega)$, where $(\phi,\psi)$ is an eigenvector of the pencil ${\cal L}$ corresponding to the eigenvalue
$\lambda_{-2}=-1-\lambda_2$. Obviously $(u,p)\in E_\beta^2(K)\times V_\beta^1(K)$ since $\beta>\mbox{Re}\, \lambda_2 + \frac 32$. Furthermore,
$(s-\Delta)u + \nabla p \in E_\gamma^2(K)$ and $\nabla\cdot u \in X_\gamma^1(K)$ if $\gamma> \mbox{Re}\, \lambda_2 -\frac 12$.
Suppose that $\max(\frac 52 \, ,\mbox{Re}\lambda_2 -\frac 12) < \gamma <\mbox{Re}\, \lambda_2 + \frac 32$. Then there exists a vector function $(v,q) \in E_\gamma^2(K)\times
V_\gamma^1(K)$ satisfying (\ref{2Bl4}). Since $\mu_{-2}=-1-\mu_2 <-\beta-\frac 12 < -\gamma-\frac 12 < \mu_{-1}=-1$, it follows that
$(v,q)\in E_\beta^2(K)\times V_\beta^1(K)$. Obviously, $u \not \in E_\gamma^2(K)$. Hence, $(u-v,p-q)$ is a nonzero element of
$\mbox{ker}\, A_\beta$. The proof is complete. \hfill $\Box$ \\

Finally, we prove the following regularity assertion for the solutions of the problem (\ref{par3}) which improves
Lemma \ref{l1}.

\begin{Le} \label{l2}
Suppose that $\mbox{\em Re}\, s\ge 0$, $s\not=0$, $\lambda_1=1$ and that $\lambda_1$ is a simple eigenvalue of the pencil ${\cal L}(\lambda)$.
Furthermore, we assume that $(u,p) \in E_\beta^2(K)\times V_\beta^1(K)$ is a solution of the problem {\em (\ref{par3})} with the data
$f\in E_\beta^0(K)\cap E_\gamma^0(K)$, $g\in X_\beta^1(K)\cap X_\gamma^1(K)$, where
\[
\frac 12 -\min( \mbox{\em Re}\, \lambda_2,\mu_2+1) < \beta,\gamma <\frac 12 +\min (\mu_2, \mbox{\em Re} \lambda_2+1),\quad \beta,\gamma
  \not\in \Big\{-\frac 12\, , \, \frac 12 \, ,\, \frac 52\Big\}.
\]
In the case $\max(\beta,\gamma)> \frac 12$, $\min(\beta,\gamma)< \frac 52 $, we assume in addition that $g$ satisfies the condition {\em (\ref{tildeX})}. Then
$u\in E_\gamma^2(K)$ and $p-c\in V_\gamma^1(K)$ with some constant $c=c(s)$. Here, $c=0$ if $(\beta+\frac 12)\, (\gamma+\frac 12)>0$.
\end{Le}

P r o o f.
First note that for $\max(\beta,\gamma)> \frac 12$, $\min(\beta,\gamma)< \frac 52$, there exist a number $\beta'$ between $\beta$ and $\gamma$ such
that $\frac 12 < \beta'<\frac 52$. Since $X_\beta^1(K)\cap X_\gamma^1(K) \subset X_{\beta'}^1(K)$, the integral of $g$ over $K$ exists in this case.
By Theorem \ref{Bt1}, there exists a solution $(v,q) \in E_\gamma^2(K)\times V_\gamma^1(K)$ of the problem (\ref{par3}). We show that $u=v$ and
$p-q$ is constant, $p-q=0$ if $(\beta-\frac 12)\, (\gamma-\frac 12)>0$. Without loss of generality, we may assume that $\beta>\gamma$.
We consider the following cases.

1) $-\mu_2-\frac 12 < \gamma<\beta<-\frac 12$ or $-\frac 12 < \gamma < \beta < \frac 12$
or $\frac 12 < \gamma<\beta < \mu_2+\frac 12$. In this case, the interval $-\beta-\frac 12 \le \lambda \le -\gamma-\frac 12$ does not contain eigenvalues
of the pencil ${\cal N}(\lambda)$ and Lemma \ref{l5b} implies $v\in E_\beta^2(K)$, $q\in V_\beta^1(K)$. Since the solution in $E_\beta^2(K)\times V_\beta^1(K)$
is unique it follows that $v=u$ and $q=p$.

2) $\frac 12 - \mbox{Re}\, \lambda_2 < \gamma<\beta < \frac 52$ or $\frac 52 < \gamma < \beta < \frac 32 + \mbox{Re}\, \lambda_2$.
Then the strip $\frac 12-\beta \le \mbox{Re}\, \lambda \le \frac 12-\gamma$ contains at most the simple eigenvalue $\lambda_1=1$ of the
pencil ${\cal L}(\lambda)$ with the constant eigenvector $(0,1)$. If $\gamma\ge \beta-2$, then
\[
-\Delta u + \nabla p = f-su \in V_\gamma^0(K), \quad -\nabla\cdot u = g \in V_\gamma^1(K).
\]
Using regularity results for solutions of elliptic problems in the spaces $V_\beta^l(K)$ (see, e.~g., \cite[Chapter 3, Theorem 5.5]{np-94}),
we conclude that $u \in V_\gamma^2(K)$ and $p-c \in V_\gamma^1(K)$ with some constant $c$. If $\gamma<\beta<-\frac 12$ or $\beta>\gamma>-1/2$, then
even $p\in V_\gamma^1(K)$, since $\lambda_1=1$ lies outside the strip $\frac 12-\beta \le \mbox{Re}\, \lambda \le \frac 12-\gamma$.
Obviously, $E_\beta^2(K) \cap V_\gamma^2(K) \subset E_\gamma^2(K)$ for $\gamma<\beta$. Thus, it is shown that $u\in E_\gamma^2(K)$ and $p-c
\in V_\gamma^1(K)$ if $\gamma \ge \beta-2$. Repeating this argument, we obtain the same result for $\gamma<\beta-2$. Since the solution is unique
in the space $E_\gamma^2(K) \times V_\gamma^1(K)$, it follows that $u=v$ and $p-c=q$. If $\gamma<\beta<-\frac 12$ or $\beta>\gamma>-1/2$, then
even $p=q$.

3) $\mu_2>2$ and $\frac 12 -\min( \mbox{Re}\, \lambda_2,\mu_2+1) < \gamma < \frac 52 < \beta < \frac 12 +\min (\mu_2, \mbox{Re} \lambda_2+1)$.
Then we can choose a number $\gamma'$ such that $\max(\gamma,\frac 12) < \gamma' < \frac 52$. Since $f\in E_{\gamma'}^0(K)$ and $g\in X_{\gamma'}^1(K)$,
there exists a solution $(v',q')\in E_{\gamma'}^2(K)\times V_{\gamma'}^1(K)$ of the problem (\ref{par3}). As was shown in part 1), we obtain
$(u,p)=(v',q)$. Hence by 2), $(u,p) = (v,q)$ if $\gamma > -\frac 12$ and $(u,p)=(v,q-c)$ if $\gamma < -\frac 12$. The proof of the lemma is complete.
\hfill $\Box$

\section{Behavior of solutions of the parameter-depending problem at infinity}

Suppose that $f\in E_\beta^0(K)$ and $g\in X_\beta^1(K)$, where $\frac 12-\lambda_1 < \beta < \frac 12$. Then by Theorem \ref{t3}
there exists a unique solution $(u,p) \in E_\beta^2(K) \times V_\beta^1(K)$ of the problem (\ref{par3}).
By Lemma \ref{l1}, this solution belongs to the space $E_\gamma^2(K) \times V_\gamma^1(K)$ if $f\in E_\gamma^0(K)$, $g\in X_\gamma^1(K)$
and $\beta<\gamma<\frac 12$. However, this is not true in general if $\gamma>\frac 12$. We show that then the solution is a sum
of some singular terms and a remainder $(v,q) \in E_\gamma^2(K)\times V_\gamma^1(K)$.

\subsection{Special solutions of the parameter-depending problem}

In the sequel, let $\nu(x)$ denote the distance of the point $x$ from the boundary $\partial K$. Obviously, the function $\nu$
is positively homogeneous of degree 1. In the neighborhood $\nu(x)<\delta|x|$ of the boundary $\partial K$ with sufficiently small $\delta$,
the function $\nu$ is two times continuously differentiable and satisfies the equality $|\nabla \nu|=1$. Furthermore,
the vector $\nabla \nu(x)$ is orthogonal to $\partial K$ at any point $x\in \partial K\backslash\{ 0\}$.  For an arbitrary vector function $v$ in the neighborhood
$\nu(x)<\delta |x|$ of $\partial K$, we define
\[
v_\nu = v\cdot \nabla\nu \quad\mbox{and}\quad v_\tau = v-v_\nu\, \nabla \nu.
\]
Obviously $v_\tau\cdot\nabla\nu=0$ near $\partial K$.

Now, let $\mu$ be an eigenvalue of the pencil ${\cal N}(\lambda)$, and let $\phi$ be an eigenfunction corresponding to this eigenvalue.
Then the function
\begin{equation} \label{p0}
p_0(x)=r^\mu\phi(\omega)
\end{equation}
is a solution of the Neumann problem $\Delta p_0=0$ in $K$, $\nabla p_0\cdot \nabla \nu =0$ on $\partial K\backslash \{ 0\}$.
Furthermore, let $\chi$ be a two times continuously differentiable function on $(0,\infty)$
such that $\chi(r)=1$ for $2r<\delta$ and $\chi(r)=0$ for $r>\delta$. We define
\begin{equation} \label{u0}
u_0(x,s) = s^{-1}\, \Big( v^{(0)}(x) - \chi\Big( \frac \nu r \Big) \, e^{-\nu\sqrt{s}}\, v^{(0)}_\tau(x)\Big),\quad \mbox{where }\ v^{(0)}= - \nabla p_0
\end{equation}
and $\sqrt{s}$ is the square root of $s$ with positive real part. Then the following assertion holds.

\begin{Le}  \label{Cl1}
Let $p_0$ and $u_0$ be the functions {\em (\ref{p0})} and {\em (\ref{u0})}, respectively. Then
\[
(s-\Delta)\, u_0 + \nabla p_0 =  \chi\Big( \frac\nu r \Big)\, e^{-\nu\sqrt{s}}\, s^{-1/2}\,  \Big(
  - v^{(0)}_\tau \Delta\nu- 2\sum_{j=1}^3 \frac{\partial\nu}{\partial x_j}\, \frac{\partial v^{(0)}_\tau}{\partial x_j}
  + s^{-1/2}\, \Delta\, v^{(0)}_\tau\Big) + R_1
\]
and
\[
\nabla\cdot u_0 = -  \chi\Big( \frac\nu r \Big)\, e^{-\nu\sqrt{s}}\, s^{-1}\, \nabla\cdot v^{(0)}_\tau + R_2,
\]
where  $R_1,R_2$ are continuous in $K$ and satisfy the estimates
\begin{equation} \label{1Cl1}
|R_1| \le c\, s^{-1}\, r^{\mu-3}\, e^{-\delta r \, \mbox{\em \scriptsize Re} \sqrt{s}/3},\quad
|R_2| \le c\, s^{-1}\, r^{\mu-2}\, e^{-\delta r \, \mbox{\em \scriptsize Re} \sqrt{s}/3}
\end{equation}
with a constant $c$ independent of $s$ and $r$. Furthermore, $u_0=0$ on $\partial K \backslash\{0\}$.
\end{Le}

The remainders $R_1$ and $R_2$ in Lemma \ref{Cl1} are
\[
R_1= s^{-1}\, \big[ \Delta,\chi(\frac\nu r )\big]\,  e^{-\nu\sqrt{s}}\,  v^{(0)}_\tau\quad\mbox{and}\quad
R_2= -s^{-1}\, e^{-\nu\sqrt{s}}\, v^{(0)}_\tau\cdot \nabla \chi(\frac\nu r ).
\]
Here $[\Delta,\chi]$ denotes the commutator of $\Delta$ and $\chi$. Since $\nabla \chi(\frac\nu r )$ is equal to zero for $2\nu<\delta r$, the terms $R_1,R_2$ satisfy the estimate (\ref{1Cl1}). Obviously, the terms $v^{(0)}_\tau\Delta\nu$,
$\frac{\partial\nu}{\partial x_j}\, \frac{\partial v^{(0)}_\tau}{\partial x_j}$ and $\nabla\cdot v^{(0)}_\tau$ are positively homogeneous of degree $\mu-2$,
while $\Delta\, v^{(0)}_\tau$ is positively homogeneous of degree $\mu-3$. \\

It follows from the last lemma that
\[
\big| (s-\Delta)\, u_0 + \nabla p_0 \big| + |s|^{1/2}\, \big| \nabla\cdot u_0\big| \le c\, |s|^{-1/2}\, r^{\mu-2}\ \mbox{ for }r>|s|^{-1/2}.
\]
Next, we construct functions $U_N$ and $P_N$ with the leading terms $u_0$ and $p_0$, respectively, such that
\[
\big| (s-\Delta)\, U_N + \nabla P_N \big| + |s|^{1/2}\, \big| \nabla\cdot U_N\big| \le c\, |s|^{-(N+1)/2}\, r^{\mu-N-2}\,
  (1+|\log r|^k)\ \mbox{ for }r>|s|^{-1/2},
\]
where $k\le N$. We introduce the polynomials
\[
P_1(\nu)=-k! \sum_{j=0}^k \frac{\nu^j}{j!}\, ,\quad P_2(\nu)= 2^{-k-2}\, k! \sum_{j=1}^{k+1} \frac{(2\nu)^j}{j!}\, ,\quad
P_3(\nu)=-k! \sum_{j=0}^{k-1} (k-j+1)\, \frac{\nu^j}{j!}
\]
for integer $k\ge 0$ ($P_3=0$ if $k=0$) which satisfy the equalities
\[
P'_1(\nu)-P_1(\nu)=\nu^k, \quad P''_2(\nu)-2P'_2(\nu)=-\nu^k,\quad P'_3(\nu)-P_3(\nu)=P''_1(\nu)-2\, P'_1(\nu).
\]
Then the following lemma holds.

\begin{Le} \label{Cl2}
Let $f$ and $g$ be positively homogeneous of degree $\mu$ in the neighborhood $\nu(x)<\delta|x|$ of $\partial K$.
Then the functions
\[
u = s^{-1/2}\, e^{-\nu\sqrt{s}}\big( P_1(\nu\sqrt{s})\, g\, \nabla\nu + P_2(\nu\sqrt{s})\, f_\tau\big), \quad
p = e^{-\nu\sqrt{s}}\big(  P_3(\nu\sqrt{s})\, g + P_1(\nu\sqrt{s})\, f_\nu \big)
\]
satisfy the equations
\[
(s-\Delta)\, u + \nabla p = \sqrt{s}\, e^{-\nu\sqrt{s}}\, \big( (\nu\sqrt{s})^k f + R_1\big), \quad
\nabla\cdot u = e^{-\nu\sqrt{s}}\, \big( (\nu\sqrt{s})^k g + R_2\big)
\]
in the neighborhood $\nu(x)<\delta|x|$ of $\partial K$, where $R_1,R_2$ have the form
\begin{eqnarray*}
R_1 & = & s^{-1/2}\, r^{\mu-1}\, \sum_{j=0}^{k+1} (\nu\sqrt{s})^j\, a_j(\omega) + s^{-1}\, r^{\mu-2}\, \sum_{j=1}^{k+1} (\nu\sqrt{s})^j\, b_j(\omega),\\
R_2 & = & s^{-1/2}\, r^{\mu-1}\, \sum_{j=0}^{k+1} (\nu\sqrt{s})^j\, c_j(\omega)
\end{eqnarray*}
Furthermore, $u_\tau=0$ and $u_\nu=s^{-1/2}\, P_1(0)\, g$ on $\partial K\backslash\{ 0\}$.
\end{Le}

P r o o f.
Since $P_2(0)=0$, we have $u_\tau=0$ and $u_\nu=s^{-1/2}\, P_1(0)\, g$ on $\partial K\backslash\{ 0\}$.
One easily checks that
\begin{eqnarray*}
\nabla\cdot u & = & - e^{-\nu\sqrt{s}} \, P_1(\nu\sqrt{s})\, g
  + s^{-1/2}\, e^{-\nu\sqrt{s}}\, \nabla \cdot \big( P_1(\nu\sqrt{s})\,  g\, \nabla\nu + P_2(\nu\sqrt{s})\, f_\tau\big) \\
& = & e^{-\nu\sqrt{s}} \, \Big( (\nu\sqrt{s})^k \, g + s^{-1/2}\, P_1(\nu\sqrt{s})\, \nabla\cdot\big(g\, \nabla\nu\big)
  + s^{-1/2}\, P_2(\nu\sqrt{s})\nabla\cdot f_\tau\Big).
\end{eqnarray*}
Here, the functions $\nabla\cdot\big(g\, \nabla\nu\big)$ and $\nabla\cdot f_\tau$ are homogenous of degree $\mu-1$. Furthermore,
\[
\nabla p = \sqrt{s}\, e^{-\nu\sqrt{s}}\Big( (\nu\sqrt{s})^k  f_\nu\, \nabla\nu + \big( P'_3(\nu\sqrt{s})-P_3(\nu\sqrt{s})\big)\, g\nabla\, \nu + R\Big)
\]
where $R= s^{-1/2} \, P_1(\nu\sqrt{s})\, \nabla f_\nu + s^{-1/2}\, P_3(\nu\sqrt{s})\, \nabla g$.
Using the equalities $(s-\Delta)\, e^{-\nu\sqrt{s}} = \sqrt{s}\, e^{-\nu\sqrt{s}}\, \Delta\nu$ and
$\Delta q(\nu)=q''(\nu)+q'(\nu)\, \Delta \nu$, we get
\begin{eqnarray*}
(s-\Delta)\, u &  = & e^{-\nu\sqrt{s}}\, \Delta\nu \,\big( P_1(\nu\sqrt{s})\, g\nabla\nu+ P_2(\nu\sqrt{s})\, f_\tau\big)\, \\
&&   + 2\, e^{-\nu\sqrt{s}} \sum_{j=1}^3 \frac{\partial \nu}{\partial x_j}\,
  \frac{\partial}{\partial x_j} \Big( P_1(\nu\sqrt{s})\, g\nabla\nu+ P_2(\nu\sqrt{s})\, f_\tau\big) \\
&& -s^{-1/2}\, e^{-\nu\sqrt{s}}\, \Delta\big( P_1(\nu\sqrt{s})\, g\nabla\nu+ P_2(\nu\sqrt{s})\, f_\tau\big)\\
& = & \sqrt{s}\, e^{-\nu\sqrt{s}}\, \Big( (\nu\sqrt{s})^k\, f_\tau + \big( 2P'_1(\nu\sqrt{s})-P''_1(\nu\sqrt{s})\big)\, g\nabla \nu +R'\Big),
\end{eqnarray*}
where
\begin{eqnarray*}
R' & = & s^{-1/2}\, \big(  P_2(\nu\sqrt{s})-P'_2(\nu\sqrt{s})\big)\,  \Big( f_\tau\, \Delta\nu + 2\sum_{j=1}^3 \frac{\partial\nu}{\partial x_j} \,
  \frac{\partial f_\tau}{\partial x_j}\Big) - s^{-1}\, P_2(\nu\sqrt{s})\, \Delta f_\tau  \\
&& - s^{-1/2}\, (\nu\sqrt{s}))^k\, \Big( g\Delta\nu\, \nabla\nu + 2\sum_{j=1}^3 \frac{\partial\nu}{\partial x_j} \,
  \frac{\partial(g\nabla\nu)}{\partial x_j}\Big)  - s^{-1}P_1(\nu\sqrt{s})\, \Delta(g\nabla\nu) \, .
\end{eqnarray*}
Consequently,
\[
(s-\Delta)\, u + \nabla p = \sqrt{s}\, e^{-\nu\sqrt{s}}\, \big( (\nu\sqrt{s})^k f + R + R'\big).
\]
This proves the lemma. \hfill $\Box$ \\

Next, we construct special singular functions with the leading part $(u_0,p_0)$ defined by
(\ref{p0}) and (\ref{u0}).

\begin{Le} \label{Cl3}
Let $\mu$ be an eigenvalue of the pencil ${\cal N}(\lambda)$ with the eigenfunction $\phi$, and let
$p_0$ and $u_0$ be the functions {\em (\ref{p0})} and {\em (\ref{u0})}, respectively.
There exist functions $u_j,p_j$ of the form
\[
u_j(x,s) = s^{-(j+2)/2}\, r^{\mu-j-1} \sum_{k=0}^j (\log r)^k \, \Phi_{j,k}(\omega), \quad
p_j(x,s) = s^{-j/2}\, r^{\mu-j} \sum_{k=0}^j (\log r)^k \, \Psi_{j,k}(\omega),
\]
$j=1,2,\ldots$, and functions $w_j,q_j$ of the form
\begin{eqnarray*}
w_j(x,s) &= & s^{-(j+2)/2}\, r^{\mu-j-1} e^{-\nu\sqrt{s}}\, \sum_{i=0}^j \sum_{k=0}^j (\nu\sqrt{s})^i\, (\log r)^k \, \Phi_{i,j,k}(\omega), \\
q_j(x,s) &= & s^{-(j+1)/2}\, r^{\mu-j-1} e^{-\nu\sqrt{s}}\, \sum_{i=0}^j \sum_{k=0}^j (\nu\sqrt{s})^i\, (\log r)^k \, \Psi_{i,j,k}(\omega)
\end{eqnarray*}
such that the functions
\[
U_N = \sum_{j=0}^N u_j - \chi\Big(\frac{\nu}r \Big)\, \sum_{j=1}^{N} w_j \, , \quad
P_N = \sum_{j=0}^N p_j - \chi\Big(\frac{\nu}r \Big)\,  \sum_{j=1}^{N} q_j
\]
have the following properties for $N=0,1,2,\ldots$.

{\em 1)} There are the representations
\[
(s-\Delta)\, U_N + \nabla P_N =  \chi\Big(\frac{\nu}r \Big)\, f + R_1, \quad \nabla\cdot U_N = \chi\Big(\frac{\nu}r \Big)\, g + R_2,
\]
where $f$ and $g$ are finite sums of terms of the form
\begin{equation} \label{1Cl3}
s^{(1-d)/2}\, r^{\mu-d}\,  e^{-\nu\sqrt{s}}\, (\nu\sqrt{s})^j\, (\log r)^k\, \Phi(\omega),\quad
\end{equation}
and
\begin{equation} \label{2Cl3}
s^{-d/2}\, r^{\mu-d}\, e^{-\nu\sqrt{s}}\, (\nu\sqrt{s})^j\, (\log r)^k\, \Psi(\omega),
\end{equation}
respectively, with nonnegative integers $d\ge N+2$ and $j,k\le N$. The remainders $R_1,R_2$ vanish for $2\nu<\delta r$ and for $\nu>\delta r$ and
satisfy the estimate
\begin{equation} \label{3Cl3}
|R_1| + |s|^{1/2}\, |R_2| \le c\, |s|^{(1-\mu)/2}\, \big( 1+\big|\log |s|\big|^N\big)\, e^{-\delta r \, \mbox{\em \scriptsize Re} \sqrt{s}/3}
\end{equation}
for $r>|s|^{-1/2}$.

{\em 2)} The function $U_N$ satisfies the boundary condition $U_N=0$ on $\partial K \backslash \{ 0\}$.
\end{Le}

P r o o f.
For the case $N=0$, we refer to Lemma \ref{Cl1}. We assume that we already constructed the functions $U_N$, $P_N$ for a certain $N\ge 0$.
Then
\[
(s-\Delta)\, U_N + \nabla P_N =  \chi\Big( \frac\nu r \Big)\, (f_1+f_2) + R_1, \quad
\nabla\cdot U_N = \chi\Big( \frac\nu r \Big)\, (g_1+g_2) + R_2.
\]
Here $f_1,g_1$ are finite sums of terms of the form (\ref{1Cl3}), (\ref{2Cl3}), respectively, with $d=N+2$ and integers $j,k\le N$, while
$f_2,g_2$ are functions of the same form with $d\ge N+3$ and $j,k\le N$. The remainders $R_1,R_2$ satisfy (\ref{3Cl3}).

By Lemma \ref{Cl2}, there exist functions $v,q$ of the form
\begin{eqnarray*}
v(x,s) & = & s^{-(N+3)/2}\, r^{\mu-N-2}\, e^{-\nu\sqrt{s}} \sum_{j=0}^{N+1} \sum_{k=0}^N (\nu\sqrt{s})^j \, (\log r)^k\, \Phi_{j,k}(\omega), \\
q(x,s) & = & s^{-(N+2)/2}\, r^{\mu-N-2}\, e^{-\nu\sqrt{s}} \sum_{j=0}^{N+1} \sum_{k=0}^N (\nu\sqrt{s})^j \, (\log r)^k\, \Psi_{j,k}(\omega)
\end{eqnarray*}
in the neighborhood $\nu(x)<\delta|x|$ of $\partial K$ such that
\[
(s-\Delta)\, v + \nabla q = f_1+f_3\, , \quad \nabla\cdot v = g_1 + g_3 \, ,
\]
where $f_3$ and $g_3$ are finite sums of terms of the the form (\ref{1Cl3}), (\ref{2Cl3}), respectively, with $d\ge N+3$ and integers $j\le N+1$, $k\le N$.
Furthermore,
\[
v_\tau =0, \quad v_\nu = s^{-(N+3)/2}\, r^{\mu-N-2} \sum_{k=0}^N  \, (\log r)^k\, \Phi_{0,k}(\omega)\cdot \nabla\nu \quad
  \mbox{on }\partial K \backslash \{ 0\}.
\]
Consequently, there exists a solution $p'$ of the Neumann problem
\[
\Delta p' =0\ \mbox{ in }K,\quad \nabla p'\cdot \nabla\nu = -s\, v_\nu \ \mbox{ on }\partial K\backslash\{ 0\}
\]
which has the form
\[
p'(x,s) = s^{-(N+1)/2}\, r^{\mu-N-1} \sum_{k=0}^{N+1} (\log r)^k \, \psi_k(\omega)
\]
(see, e.~g., \cite[Lemma 6.1.13]{kmr-97}). Here $\psi_{N+1}=0$ if $\mu-N-1$ is not an eigenvalue of the pencil ${\cal N}(\lambda)$. We set
\[
v' = -  \nabla p' \quad\mbox{and}\quad u' = s^{-1}\, \Big( v' - \chi\Big( \frac{\nu}r \Big) e^{-\nu\sqrt{s}}\, v'_\tau \Big) .
\]
Analogously to Lemma \ref{Cl1}, we get the representations
\[
(s-\Delta) u' + \nabla p' = \chi\Big( \frac{\nu}r \Big)\, F + R'_1, \quad
  \nabla\cdot u'=  \chi\Big( \frac{\nu}r \Big)\, G + R'_2 \ \mbox{ in }K,
\]
where $F$, $G$ are finite sums of terms of the form (\ref{1Cl3}) and (\ref{2Cl3}), respectively, with $d\ge N+3$, $j=0$, $k\le N+1$,
and $R'_1,R'_2$ satisfy the estimate (\ref{3Cl3}). Furthermore,
\[
u'_\tau=0, \ \ u'_\nu = s^{-1}\,v'_\nu =-s^{-1}\, \nabla p'\cdot \nabla\nu = v_\nu \ \mbox{ on }\partial K\backslash \{ 0\}.
\]
We consider the functions
\[
U_{N+1} = U_N + u'- \chi\Big(\frac{\nu}r \Big)\, v = U_N + u_{N+1}- \chi\Big(\frac{\nu}r \Big)\, w_{N+1}
\]
and
\[
P_{N+1} = P_N + p'- \chi\Big(\frac{\nu}r \Big)\, q = P_N + p_{N+1}- \chi\Big(\frac{\nu}r \Big)\, q_{N+1},
\]
where $u_{N+1}= s^{-1}v'$, $w_{N+1}= v+s^{-1}\, e^{-\nu\sqrt{s}}\, v'_\tau$, $p_{N+1}=p'$ and $q_{N+1}=q$. Obviously,
$u_{N+1},p_{N+1},w_{N+1}$ and $q_{N+1}$ have the desired form, and the function $U_{N+1}$ satisfies the boundary condition
$U_{N+1}=0$ on $\partial K \backslash \{ 0\}$. Furthermore, we obtain
\[
(s-\Delta)\, U_{N+1} + \nabla P_{N+1} = \chi\Big( \frac{\nu}r \Big)\, (f_2-f_3+F) +R_1+R'_1 +[\Delta,\mbox{$\chi(\frac\nu r )$}]\, v
  -q\, \nabla\mbox{$\chi(\frac\nu r )$}
\]
and
\[
\nabla\cdot U_{N+1} = \chi\Big( \frac{\nu}r \Big)\, (g_2-g_3+G) +R_2+R'_2 -v\cdot \nabla\mbox{$\chi(\frac\nu r )$}.
\]
This proves the lemma. \hfill $\Box$

\begin{Rem} \label{Cr1}
{\em 1) In the case $\mu=0$, where $u_0=0$ and $p_0$ is constant, we can obviously set $U_N=u_0=0$ and $P_N=p_0$ for arbitrary $N$.

2) Logarithmic terms in the representation of $U_N$ and $P_N$ appear only in the case $N>0$ if at least one of the numbers
$\mu-1,\mu-2,\ldots,\mu-N$ is an eigenvalue of the pencil ${\cal N}(\lambda)$. In fact, the logarithm appears at most with an exponent
$N'\le N$, where $N'$ is the number of eigenvalues of the pencil ${\cal N}(\lambda)$ in the set $\{ \mu-1,\mu-2,\ldots,\mu-N\}$.}
\end{Rem}

In the following, $\eta$ is a two times continuously differentiable function on $(0,\infty)$, $\eta(r)=0$ for $r<1/2$, $\eta(r)=1$ for $r>1$.
Furthermore, we define
\[
\eta_s(x) = \eta(|s|\, r^2).
\]
Then the following assertion holds.

\begin{Le}  \label{Cl4}
Let $\mu$ be an eigenvalue of the pencil ${\cal N}(\lambda)$ with the eigenfunction $\phi$.
Furthermore, let $U_N,P_N$ be the functions described in Lemma {\em \ref{Cl3}}. Then
\begin{equation} \label{2Cl4}
\big \| (s-\Delta)\, (\eta_s U_N) + \nabla (\eta_s P_N)\big\|^2_{V_\beta^0(K)} \le c\, |s|^{-\beta-\mu-1/2}\, \big( 1+\big|\log|s|\big|^{2N'}\big)
\end{equation}
for $\beta+\mu< N+1$ and
\[
\big\| \nabla\cdot (\eta_s U_N )\big\|^2_{V_\beta^1(K)}+ |s|^2\, \big\| \nabla\cdot (\eta_s U_N )\big\|^2_{(V_{-\beta}^1(K))^*}
  \le c\, |s|^{-\beta-\mu-1/2}\, \big( 1+\big|\log|s|\big|^{2N''}\big)
\]
for $\beta+\mu<N+ \frac 12$, where $c$ is a constant independent of $s$.
Here, $N'$ is the number of eigenvalues of the pencil ${\cal N}(\lambda)$ in the set $\{ \mu-1,\mu-2,\ldots,\mu-N\}$,
and $N''\le N'+1$. If $\mu-1,\mu-2,\ldots,\mu-N$ are not eigenvalues of the pencil ${\cal N}(\lambda)$, then $N'=N''=0$.
\end{Le}

P r o o f.
By Lemma \ref{Cl3},
\[
\eta_s\,  \big( (s-\Delta)\, U_N + \nabla P_N\big) = \eta_s\, \chi(r^{-1}\nu)\, e^{-\nu\sqrt{s}}\, f+ \eta_s\, R_1
\]
and \ $\eta_s\, \nabla\cdot U_N = \eta_s\, \chi(r^{-1}\nu)\, e^{-\nu\sqrt{s}}\, g+\eta_s\, R_2$,
where $f,g$  are sums of terms of the form (\ref{1Cl3}) and (\ref{2Cl3}), respectively, and $R_1,R_2$ satisfy the estimate
(\ref{3Cl3}). Obviously,
\[
\| \eta_s \chi e^{-\nu\sqrt{s}} f\|^2_{V_\beta^0(K)} \le c\, |s|^{-N-1}\int_K |\eta_s|^2 \, \Big| \chi\Big( \frac \nu r\Big)\Big|^2 \,
  r^{2(\beta+\mu-N-2)}\, \big( 1+|\log r|^{2N'}\big)\,  e^{-\nu\, \mbox{\scriptsize Re}\, \sqrt{s}}\, dx.
\]
Here, $N'$ is the number of eigenvalues of the pencil ${\cal N}(\lambda)$ in the set $\{ \mu-1,\mu-2,\ldots,\mu-N\}$ (cf. Remark \ref{Cr1}). Since
\begin{equation} \label{ineq}
\int_\Omega \Big| \chi\Big( \frac{\nu(x)}r \Big)\Big|^2 \,   e^{-\nu(x)\, \mbox{\scriptsize Re}\, \sqrt{s}}\, d\omega
  =  \int_\Omega \big| \chi\big( \nu(\omega) \big)\big|^2 e^{-r\, \nu(\omega)\, \mbox{\scriptsize Re}\, \sqrt{s}}\, d\omega \le
  c\, |s|^{-1/2}\, r^{-1},
\end{equation}
we get
\begin{eqnarray*}
\| \eta_s \chi e^{-\nu\sqrt{s}}\, f\|^2_{V_\beta^0(K)} & \le & c\, |s|^{-N-3/2}\, \int_{(2|s|)^{-1/2}}^\infty (1+|\log r|^{2N'})\, r^{2\beta+2\mu-2N-3}\, dr \\
& \le & c\, |s|^{-\beta-\mu-1/2}\, \big( 1+\big|\log|s|\big|^{2N'}\big) \quad\mbox{if } \beta+\mu-N<1.
\end{eqnarray*}
The same estimate holds for $\eta_s R_1$. Thus, (\ref{2Cl4}) holds.

Analogously, we obtain
\[
\|  \nabla\cdot (\eta_s U_N)\|^2_{(V_{-\beta}^1(K))^*} \le \| \nabla\cdot (\eta_s U_N)\|^2_{V_{\beta+1}^0(K)}
  \le c\, |s|^{-\beta-\mu-5/2}\, \big( 1+\big|\log|s|\big|^{2N'}\big)
\]
and
\[
\| \nabla\cdot (\eta_s U_N)\|^2_{V_\beta^1(K)}  \le c\, |s|^{-\beta-\mu-1/2}\, \big( 1+\big|\log|s|\big|^{2N'}\big)
\]
if $\beta+\mu<N$. Suppose that $N\le \beta+\mu<N+\frac 12$. Then
\[
\| \nabla\cdot (\eta_s U_{N+1})\|^2_{V_\beta^1(K)} + |s|^2\, \| \nabla\cdot (\eta_s U_{N+1})\|^2_{(V_{-\beta}^1(K))^*}
  \le c\, |s|^{-\beta-\mu-1/2}\, \big( 1+\big|\log|s|\big|^{2N''}\big),
\]
where $N''$ is the number of eigenvalues of the pencil ${\cal N}(\lambda)$ in the set $\{ \mu-1,\mu-2,\ldots,\mu-N-1\}$.
Furthermore, one can easily show that
\[
\eta_s (U_{N+1}-U_N) = \eta_s \big( u_{N+1}-\chi(r^{-1}\nu)\, w_{N+1}\big) \in E_\beta^2(K) \cap \stackrel{\circ}{E}\!{}_\beta^1(K)
\]
and
\[
\big\| \eta_s (U_{N+1}-U_N)\big\|^2_{V_\beta^2(K)} + |s|^2 \, \big\| \eta_s (U_{N+1}-U_N)\big\|^2_{V_\beta^0(K)}
  \le c\, |s|^{-\beta-\mu-1/2}\, \big( 1+\big|\log|s|\big|^{2N''}\big).
\]
This implies
\begin{eqnarray*}
&& \| \nabla\cdot (\eta_s U_{N+1}-\eta_s U_N)\|^2_{V_\beta^1(K)} + |s|^2\, \|  \nabla\cdot (\eta_s U_{N+1}-\eta_s U_N)\|^2_{(V_{-\beta}^1(K))^*}\\
&&  \le c\, |s|^{-\beta-\mu-1/2}\, \big( 1+\big|\log|s|\big|^{2N''}\big).
\end{eqnarray*}
Thus, the desired estimate for $\nabla\cdot (\eta_s U_N)$ holds in the case $\beta+\mu<N+\frac 12$.
If $N'=0$, i.~e., $\mu-1,\mu-2,\ldots,\mu-N$ are not eigenvalues of the pencil ${\cal N}(\lambda)$, then $U_N$ and $P_N$ do not contain logarithmic
terms (see Remark \ref{Cr1}) and we have
\[
U_N(x,s) = |s|^{-(\mu+1)/2}\, U_N\big( |s|^{1/2}x, |s|^{-1}s\big).
\]
We define $\hat{U}_{N+1}(x,s) = |s|^{-(\mu+1)/2} \, U_{N+1}\big( |s|^{1/2}x, |s|^{-1}s\big)$, i.~e.,  $\hat{U}_{N+1}(x,s)$ arises
if we replace $\log r$ by $\log (|s|^{1/2}r)$ in the representation of $U_{N+1}(x,s)$. Then
\[
\| \nabla\cdot (\eta_s \hat{U}_{N+1})\|^2_{V_\beta^1(K)} + |s|^2\, \| \nabla\cdot (\eta_s \hat{U}_{N+1})\|^2_{(V_{-\beta}^1(K))^*}
  \le c\, |s|^{-\beta-\mu-1/2}.
\]
Furthermore,
\[
\hat{U}_{N+1} - U_N = |s|^{-(\mu+1)/2} \, \Big( u_{N+1}\big( |s|^{1/2}x, |s|^{-1}s\big)- \chi\Big(\frac{\nu}r \Big)\,
  w_{N+1}\big( |s|^{1/2}x, |s|^{-1}s\big)\Big).
\]
One can easily show that
\[
\big\| \eta_s (\hat{U}_{N+1}-U_N)\big\|^2_{V_\beta^2(K)} + |s|^2 \, \big\| \eta_s (\hat{U}_{N+1}-U_N)\big\|^2_{V_\beta^0(K)}
  \le c\, |s|^{-\beta-\mu-1/2}
\]
and, consequently,
\begin{eqnarray*}
\| \nabla\cdot (\eta_s \hat{U}_{N+1}-\eta_s U_N)\|^2_{V_\beta^1(K)} + |s|^2\, \|  \nabla\cdot (\eta_s U_{N+1}-\eta_s U_N)\|^2_{(V_{-\beta}^1(K))^*}
 \le c\, |s|^{-\beta-\mu-1/2}.
\end{eqnarray*}
This proves the inequality
\[
\big\| \nabla\cdot (\eta_s U_N )\big\|^2_{V_\beta^1(K)}+ |s|^2\, \big\| \nabla\cdot (\eta_s U_N )\big\|^2_{(V_{-\beta}^1(K))^*}
  \le c\, |s|^{-\beta-\mu-1/2}
\]
for the case $N'=0$. The proof of the lemma is complete. \hfill $\Box$ \\

Using Lemma \ref{Cl4}, we can construct special solutions $(V,Q)$ of the problem
\begin{equation} \label{1Cl4}
(s-\Delta)\, V + \nabla Q =0, \quad \nabla\cdot V =0\ \mbox{ in }K, \quad V=0\ \mbox{ on } \partial K,
\end{equation}
with the leading term $\eta_s (u_0,p_0)$, where $p_0,u_0$ are the functions (\ref{p0}), (\ref{u0}) with a nonnegative eigenvalue
$\mu$ of the pencil ${\cal N}(\lambda)$.

\begin{Co} \label{Cc1}
Let $\mu$ be a nonnegative eigenvalue of the pencil ${\cal N}(\lambda)$ with the eigenfunction $\phi$, and let $N$ be the smallest integer
greater than $\mu-\lambda_1$. Furthermore, let $U_N,P_N$ are the functions described in Lemma {\em \ref{Cl3}}.
Then there exists a nontrivial solution $(V,Q)$ of the problem {\em (\ref{1Cl4})} which has the form
\[
V = \eta_s\, U_N + v, \quad  Q= \eta_s\, P_N + q,
\]
where $(v,q)\in E_\beta^2(K)\times V_\beta^1(K)$ with some $\beta\in \big(\frac 12 -\lambda_1,\, \frac 12\big)$.
\end{Co}

P r o o f.
By Lemma \ref{Cl4}, we have $(s-\Delta)(\eta_s U_N) + \nabla(\eta_s P_N) \in E_\beta^0(K)$ and $\nabla\cdot(\eta_s U_N) \in X_\beta^1(K)$
for $\beta+\mu<N+\frac 12$. Moreover, $\eta_S U_N=0$ on $\partial K \backslash\{ 0\}$. Suppose that
\[
\frac 12 - \lambda_1 < \beta < \min(N-\mu,0)+\frac 12 \, .
\]
Then by Theorem \ref{t3}, there exists a pair $(v,q)\in E_\beta^2(K)\times V_\beta^1(K)$ such that
\[
(s-\Delta)\, (\eta_s U_N+v) + \nabla(\eta P_N+q) =0, \ \ \nabla\cdot (\eta_s U_N+v)= \ \mbox{ in }K
\]
and $v=0$ on $\partial K\backslash \{ 0\}$.
Obviously, $\eta_s P_N \not\in V_\beta^1(K)$ since $\beta>-\frac 12$ and $\mu\ge 0$. Hence $\eta_s P_N+q\not=0$. This proves the corollary. \hfill $\Box$ \\

In the case $\mu\not=\mu_1=0$, the vector function $(V,Q)$ in Corollary \ref{Cc1} depends on $s$.
In the case $\mu=0$, we have $N=0$, $U_0=u_0=0$ and $P_0=p_0$ is a constant. If we choose $v=0$, $q=(1-\eta_s)\, p_0$, we get the solution
$(V,Q)=(0,p_0)$ of the problem (\ref{1Cl4}). Obviously, $q\in V_\beta^1(K)$ for arbitrary $\beta>-\frac 12$.

\subsection{Asymptotics of the solution}

Let $\mu_j$, $j=1,2,\ldots$, denote the nonnegative eigenvalues of the pencil ${\cal N}$ and let
$\phi_{j,k}$, $k=1,\ldots,\sigma_j$, be orthonormalized eigenfunctions corresponding to $\mu_j$. Then the
functions $\phi_{-j,k}=\phi_{j,k}$ are also eigenfunctions corresponding to the negative eigenvalues $\mu_{-j}=-1-\mu_j$.
For arbitrary integer $j\not=0$, and $k=1,\ldots,\sigma_j$, we set (cf. (\ref{p0}), (\ref{u0}))
\[
p_0^{(j,k)}(x)=r^{\mu_j}\, \phi_{j,k}(\omega), \quad
  u_0^{(j,k)}(x,s) = s^{-1}\, \Big( v^{(j,k)}(x) - \chi\Big( \frac \nu r \Big) \, e^{-\nu\sqrt{s}}\, v^{(j,k)}_\tau(x)\Big),
\]
where $v^{(j,k)}= - \nabla p_0^{(j,k)}$.
Furthermore, let $u_n^{(j,k)}$, $p_n^{(j,k)}$, $w_n^{(j,k)}$, and $q_n^{(j,k)}$, $n=1,2,\ldots,$ be the functions described in Lemma \ref{Cl3}
for the eigenvalue $\mu=\mu_j$ and the eigenfunction $\phi=\phi_{j,k}$. For arbitrary integer $N\ge 0$ we define (cf. Lemma \ref{Cl3})
\[
U_N^{(j,k)} = \sum_{n=0}^N u_n^{(j,k)} - \chi\Big(\frac{\nu}r \Big)\, \sum_{n=1}^{N} w_n^{(j,k)} \, , \quad
P_N^{(j,k)} = \sum_{n=0}^N p_n^{(j,k)} - \chi\Big(\frac{\nu}r \Big)\,  \sum_{n=1}^{N} q_n^{(j,k)}.
\]
Let $\mu_{-j}=-1-\mu_j$ be a negative eigenvalue of ${\cal N}(\lambda)$ and let $\gamma>-\mu_{-j}-\frac 12 = \mu_j+\frac 12$. Then we denote by
$M_{j,\gamma}$ the smallest integer such that $M_{j,\gamma}>\gamma-\mu_j-\frac 32$ and set
\[
U^{j,k,\gamma} = U_{M_{j,\gamma}}^{(-j,k)}, \quad P^{j,k,\gamma} = P_{M_{j,\gamma}}^{(-j,k)} \quad \mbox{for }j=1,2,\ldots.
\]

\begin{Le} \label{Cl4a}
Suppose that $\mbox{\em Re}\, s \ge 0$, $s\not= 0$, and that $(u,p) \in E_\beta^2(K)\times V_\beta^1(K)$ is a solution of the problem
{\em (\ref{par3})} with the data
\[
f\in E_\beta^0(K)\cap E_\gamma^0(K), \quad g\in X_\beta^1(K)\cap X_\gamma^1(K),\quad\mbox{where } \gamma>\beta>-\frac 12 \, .
\]
If the numbers $\beta-\frac 12$ and $\gamma-\frac 12$ are not eigenvalues of the pencil ${\cal N}(\lambda)$, then $(u,p)$ admits the decomposition
\[
(u,p) = \eta_s(r) \sum_{j \in J_{\beta,\gamma}} \sum_{k=1}^{\sigma_j} c_{j,k}(s)\, \big( U^{j,k,\gamma}, P^{j,k,\gamma}\big) + (v,q),
\]
where $v\in E_\gamma^2(K)$, $q\in V_\gamma^1(K)$ and $J_{\beta,\gamma}$ denotes the set of all $j$ such that $\beta-\frac 12 < \mu_j < \gamma-\frac 12$.
If $|s|=1$, then
\begin{eqnarray}
&& \| v\|_{E_\gamma^2(K)} + \| q\|_{V_\gamma^1(K)} + \sum_{j \in J_{\beta,\gamma}} \sum_{k=1}^{\sigma_j} |c_{j,k}(s)|  \nonumber \\ \label{1Cl4a}
&&  \le c\, \Big( \| f\|_{E_\gamma^0(K)} + \| g\|_{X_\gamma^1(K)} + \| u\|_{E_\beta^2(K)} + \| p\|_{V_\beta^1(K)}\Big)
\end{eqnarray}
with a constant $c$ independent of $s$.
\end{Le}

P r o o f.
We assume first that $\gamma \le \beta+\frac 12$. Then $-\frac 12 < \gamma-\mu_j-1<0$ for $j\in J_{\beta,\gamma}$ and consequently, $M_{j,\gamma}=0$.
It follows from  (\ref{par3}) that
\[
\int_K \nabla p \cdot \nabla q\, dx = \langle F,q\rangle \quad\mbox{for all }q\in V_{-\beta}^1(K),
\]
where
\[
\langle F,q\rangle = \int_K \big( (f+\Delta u)\cdot\nabla q -s\, g\, q\big)\, dx.
\]
Obviously, the functional $F$ is continuous on $V_{-\beta}^1(K)$,
\[
\| F\|_{(V_{-\beta}^1(K))^*} \le |s|\, \| u\|_{V_\beta^0(K)} + \| p\|_{V_\beta^1(K)}+  |s|\, \| g\|_{(V_{-\beta}^1(K))^*}.
\]
By \cite[Lemma 2.5]{k/r-16}, the functional $F$ is also continuous on $V_{1-\gamma}^2(K)$. For $|s|=1$, we have
\[
\| F\|_{(V_{1-\gamma}^2(K))^*} \le c\, \Big( \| f\|_{E_\gamma^0(K)} + \| g\|_{X_\gamma^1(K)} + \| u\|_{E_\beta^2(K)} \Big).
\]
Hence, it follows from \cite[Lemma 2.6]{k/r-16} that
\[
p = \sum_{j \in J_{\beta,\gamma}} \sum_{k=1}^{\sigma_j} c_{j,k}\, p_0^{(-j,k)}(x) + q',
\]
where $q'\in V_{\gamma-1}^0(K)$,
\[
\| q'\|_{V_{\gamma-1}^0(K)} + \sum_{j \in J_{\beta,\gamma}} \sum_{k=1}^{\sigma_j} |c_{j,k}| \le c\, \Big( \| F\|_{(V_{-\beta}^1(K))^*}
  + \| F\|_{(V_{1-\gamma}^2(K))^*}\Big).
\]
Since $(1-\eta_s)\, p_0^{(-j,k)} \in V_{\gamma-1}^0(K)$ for $\mu_j < \gamma-\frac 12$, this implies
\[
p(x) = \eta_s(r)  \sum_{j \in J_{\beta,\gamma}} \sum_{k=1}^{\sigma_j} c_{j,k}\, p_0^{(-j,k)}(x) +q(x)
\]
with a remainder $q\in V_{\gamma-1}^0(K)$, where
\[
\| q\|_{V_{\gamma-1}^0(K)} \le \| q'\|_{V_{\gamma-1}^0(K)} + c\, \sum_{j \in J_{\beta,\gamma}} \sum_{k=1}^{\sigma_j} |c_{j,k}|
\]
if $|s|=1$. We define
\[
v(x)= u(x)-\eta_s(r) \, \sum_{j \in J_{\beta,\gamma}} \sum_{k=1}^{\sigma_j} c_{j,k}\, u_0^{(-j,k)}(x).
\]
One easily checks that $\eta_s u_0^{(-j,k)} \in V_{\gamma-1}^1(K)$ for $j\in J_{\beta,\gamma}$. Using the imbedding
$E_\beta^2(K) \subset V_\beta^1(K)\cap V_{\beta-1}^1(K) \subset V_{\gamma-1}^1(K)$, we conclude that $v\in V_{\gamma-1}^1(K)$.
For $|s|=1$, we obtain the estimate
\[
\| v\|_{V_{\gamma-1}^1(K)} + \| q\|_{V_{\gamma-1}^0(K)}  \le  c\, \Big( \| q'\|_{V_{\gamma-1}^0(K)} + \| u\|_{E_\beta^2(K)}
  + \sum_{j \in J_{\beta,\gamma}} \sum_{k=1}^{\sigma_j} |c_{j,k}|\Big) .
\]
Furthermore, it follows from Lemma \ref{Cl4} that $(s-\Delta)\, v + \nabla q \in E_\gamma^0(K)$ and $\nabla \cdot v \in X_\gamma^1(K)$.
If $|s|=1$, then
\[
\| (s-\Delta)\, v + \nabla q \|_{E_\gamma^0(K)} + \| \nabla \cdot v \|_{X_\gamma^1(K)} \le c\,
  \big( \| f\|_{V_\gamma^0(K)} + \| g\|_{X_\gamma^1(K)} + \sum_{j \in J_{\beta,\gamma}} \sum_{k=1}^{\sigma_j} |c_{j,k}|\Big) .
\]
Moreover,  $v=0$ on $\partial K \backslash \{ 0\}$. Applying \cite[Lemma 2.4]{k/r-16}, we conclude that
$v\in E_\gamma^2(K)$ and $q\in V_\gamma^1(K)$. Furthermore, the estimate (\ref{1Cl4a}) holds if $|s|=1$.
This proves the theorem for the case $\gamma \le \beta+\frac 12$.

Suppose the assertion of the theorem is true for some $\gamma>\beta>-\frac 12$ and that
$f\in E_\beta^0(K)\cap E_{\gamma'}^0(K)$ and $g\in X_\beta^1(K)\cap X_{\gamma'}^1(K)$, where $\gamma'$ is such that
$\gamma < \gamma'\le \gamma+\frac 12$ and $\gamma'-\frac 12$ is not an eigenvalue of the pencil ${\cal N}(\lambda)$.
Since $E_\beta^0(K)\cap E_{\gamma'}^0(K) \subset E_\gamma^0(K)$ and $X_\beta^1(K)\cap X_{\gamma'}^1(K) \subset
X_\gamma^1(K)$, it follows from the induction hypothesis that
\[
(u,p) = \eta_s(r) \sum_{j \in J_{\beta,\gamma}} \sum_{k=1}^{\sigma_j} c_{j,k}\, \big( U^{j,k,\gamma}, P^{j,k,\gamma}\big) + (v',q'),
\]
where $v'\in E_\gamma^2(K)$ and $q'\in V_\gamma^1(K)$. If $|s|=1$, then
\[
\| v'\|_{E_\gamma^2(K)} + \| q'\|_{V_\gamma^1(K)} + \sum_{j \in J_{\beta,\gamma}} \sum_{k=1}^{\sigma_j} |c_{j,k}(s)|
  \le c\, \Big( \| f\|_{E_\gamma^0(K)} + \| g\|_{X_\gamma^1(K)} + \| u\|_{E_\beta^2(K)} + \| p\|_{V_\beta^1(K)}\Big).
\]
Consequently,
\begin{equation} \label{1Ct1}
(u,p) = \eta_s(r) \sum_{j \in J_{\beta,\gamma}} \sum_{k=1}^{\sigma_j} c_{j,k}\, \big( U^{j,k,\gamma'}, P^{j,k,\gamma'}\big) + (V,Q),
\end{equation}
where
\[
(V,Q) = (v',q') - \eta_s(r) \sum_{j \in J_{\beta,\gamma}} \sum_{k=1}^{\sigma_j} c_{j,k}\,
  \big( U^{j,k,\gamma'}-U^{j,k,\gamma}, P^{j,k,\gamma'}-  P^{j,k,\gamma}\big).
\]
Here, $\eta_s\, (U^{j,k,\gamma'}-U^{j,k,\gamma}) \in E_\gamma^2(K)$ since $U^{j,k,\gamma'}-U^{j,k,\gamma}$ contains only functions
$u_n^{(-j,k)}$ and $w_n^{(-j,k)}$ with index $n\ge M_{j,\gamma}+1>\gamma-\mu_j$. Analogously,
$\eta_s\, ( P^{j,k,\gamma'}-P^{j,k,\gamma}) \in V_\gamma^1(K)$. Thus, $V\in E_\gamma^2(K)$, $Q\in V_\gamma^1(K)$ and $V=0$ on
$\partial K\backslash \{ 0\}$. For $|s|=1$, we obviously get
\[
\| V \|_{E_\gamma^2(K)}+\| Q\|_{V_\gamma^1(K)} \le c\, \Big( \| f\|_{E_\gamma^0(K)} + \| g\|_{X_\gamma^1(K)} + \| u\|_{E_\beta^2(K)} + \| p\|_{V_\beta^1(K)}\Big).
\]
Let $F=(s-\Delta)\, V + \nabla Q$ and  $G=-\nabla\cdot V$. From Lemma \ref{Cl4} we conclude that $F\in V_{\gamma'}^0(K)$ and $G\in X_{\gamma'}^1(K)$.
If $|s|=1$, then
\[
\| F \|_{V_{\gamma'}^0(K)} +\| G\|_{X_{\gamma'}^1(K)} \le \| f \|_{V_{\gamma'}^0(K)} +\| g\|_{X_{\gamma'}^1(K)}
  + c\, \sum_{j \in J_{\beta,\gamma}} \sum_{k=1}^{\sigma_j} |c_{j,k}(s)|
\]
Therefore, by the first part of the proof, we have
\begin{equation} \label{2Cl4a}
(V,Q) = \eta_s(r) \sum_{j \in J_{\gamma,\gamma'}} \sum_{k=1}^{\sigma_j} c_{j,k}\, \big( u_0^{(-j,k)}, p_0^{(-j,k)}\big) + (v,q),
\end{equation}
where $v\in E_{\gamma'}^2(K)$ and $q\in V_{\gamma'}^1(K)$. For $|s|=1$, we have
\[
\| v\|_{E_{\gamma'}^2(K)} + \| q\|_{V_{\gamma'}^1(K)} + \sum_{j \in J_{\gamma,\gamma'}} \sum_{k=1}^{\sigma_j} |c_{j,k}(s)|
  \le c\, \Big( \| F\|_{V_{\gamma'}^0(K)} + \| G\|_{X_{\gamma'}^1(K)} + \| V\|_{E_\gamma^2(K)} + \| Q\|_{V_\gamma^1(K)}\Big).
\]
Combining (\ref{2Cl4a}) with (\ref{1Ct1}), we get
\[
(u,p) = \eta(r) \sum_{j \in J_{\beta,\gamma'}} \sum_{k=1}^{\sigma_j} c_{j,k}\, \big( U^{j,k,\gamma'}, P^{j,k,\gamma'}\big) + (v,q),
\]
Moreover, the desired estimate for $v,q$ and the coefficients $c_{j,k}$ holds if $|s|=1$.
Thus, the assertion of the lemma is true for all $\gamma>\beta$. \hfill $\Box$

\subsection{A formula for the coefficients}

Let $j\ge 1$, i.e., $\mu_j\ge 0$. Then we denote by $M_j$ the smallest integer greater than $\mu_j-\lambda_1$.
By Corollary \ref{Cc1}, there exist solutions $(V^{(j,k)},Q^{(j,k)})$ of the problem
(\ref{1Cl4}) which have the form
\begin{equation} \label{spec}
V^{(j,k)} = \eta_s \, U_{M_j}^{(j,k)} + v^{(j,k)}, \quad Q^{(j,k)} = \eta_s \, P_{M_j}^{(j,k)} + q^{(j,k)},
\end{equation}
where $v^{(j,k)} \in E_\beta^2(K)$ and $q^{(j,k)} \in V_\beta^1(K)$, $\frac 12 - \lambda_1 <\beta < \min(0,M_j-\mu_j)+ \frac 12$.
For the eigenvalue $\mu_1=0$, the pair $(V,Q)=(0,1)$ is a solution of the form (\ref{spec}) (see Remark \ref{Cr1}).

We introduce the bilinear forms
\[
a_s(u,p,v,q) = \int_K \Big( u\cdot \big( sv-\Delta v + \nabla q\big) - p\, \nabla\cdot v\Big)\, dx
\]
and
\[
A(u,p,v,q) = a_s(u,p,v,q)-a_s(v,q,u,p).
\]

\begin{Le} \label{Cl5}
Suppose that $(u,p) \in E_\beta^2(K) \times V_\beta^1(K)$, $(v,q)\in E_\delta^2(K)\times V_\delta^1(K)$ and $u=v=0$ on $\partial K \backslash \{ 0\}$.
If $0\le \beta+\delta\le 2$, then $A(u,p,v,q)=0$.
\end{Le}

P r o o f.
Under the assumptions of the lemma, we have $-\Delta u + \nabla p \in V_\beta^0(K)$, $-\Delta v + \nabla q \in V_\delta^0(K)$,
$\nabla\cdot u \in X_\beta^1(K)$ and $\nabla\cdot v \in X_\delta^1(K)$. Using the imbeddings
$E_\beta^2(K) \subset E_{-\delta}^0(K)$, $E_\delta^2(K)\subset V_{-\beta}^0(K)$,
\[
V_\beta^1(K) \subset V_{-\delta}^1(K)+V_{2-\delta}^1(K) \subset V_{-\delta}^1(K)+V_{1-\delta}^0(K) \subset
  V_{-\delta}^1(K)+\big( V_\delta^1(K)\big)^* = \big( X_\delta^1(K)\big)^*
\]
and $V_\delta^1(K) \subset \big( X_\beta^1(K)\big)^*$ for $0\le \beta+\delta\le 2$, we conclude that
\[
\int_K \Big( (-\Delta u + \nabla p)\cdot v - (\nabla\cdot u)\, q\Big)\, dx = \int_K \Big( u\cdot (-\Delta v+ \nabla q) -p\, \nabla\cdot v\Big)\, dx.
\]
This proves the lemma. \hfill $\Box$. \\

Furthermore, the following assertion holds.

\begin{Le} \label{Cl6}
If $i,j$ are positive integers, $\mu_j<\mu_i+1$ and $\gamma>\mu_i+1/2$, then
\[
A\big( \eta_s U_{M_j}^{(j,k)}, \eta_s P_{M_j}^{(j,k)}, \, \eta_s U^{i,l,\gamma}, \eta_s P^{i,l,\gamma}\big) = - \frac{2\mu_j+1}s \, \delta_{i,j}\, \delta_{k,l}
\]
for $k=1,\ldots,\sigma_j$, $l=1,\ldots,\sigma_i$.
\end{Le}

P r o o f.
Let $S_R$ be the intersection of the cone $K$ with the sphere $|x|=R$. Then
\begin{eqnarray*}
&& A\big( \eta_s U_{M_j}^{(j,k)}, \eta_s P_{M_j}^{(j,k)}, \, \eta_s U^{i,l,\gamma}, \eta_s P^{i,l,\gamma}\big) \\
&& = \lim_{R\to \infty} \int_{S_R} \Big( U^{i,l,\gamma}\cdot \partial_r U_{M_j}^{(j,k)} - U_{M_j}^{(j,k)}\cdot \partial_r U^{i,l,\gamma}
  + \big( P^{i,l,\gamma}\, U_{M_j}^{(j,k)} - P_{M_j}^{(j,k)}\, U^{i,l,\gamma}\big)\cdot \frac{x}R \Big)\, d\sigma .
\end{eqnarray*}
Since
\[
\big|  U^{i,l,\gamma}\cdot \partial_r U_{M_j}^{(j,k)}\big| + \big| U_{M_j}^{(j,k)}\cdot \partial_r U^{i,l,\gamma}\big| \le c R^{\mu_j-\mu_i-4}
\]
for large $R$ and $\mu_j-\mu_i < 1$, we obtain
\begin{eqnarray*}
A\big( \eta_s U_{M_j}^{(j,k)}, \eta_s P_{M_j}^{(j,k)}, \, \eta_s U^{i,l,\gamma}, \eta_s P^{i,l,\gamma}\big)
 = \lim_{R\to \infty} \int_{S_R} \big( P^{i,l,\gamma}\, U_{M_j}^{(j,k)} - P_{M_j}^{(j,k)}\, U^{i,l,\gamma}\big)\cdot \frac{x}R \, d\sigma
\end{eqnarray*}
We consider the leading terms
\[
p_0^{(-i,l)}(x) = r^{-1-\mu_i} \, \phi_{i,l}(\omega), \quad p_0^{(j,k)} = r^{\mu_j}\, \phi_{j,k}(\omega)
\]
of $P^{i,l,\gamma}$ and $P_{M_j}^{(j,k)}$ and the leading terms
\[
u_0^{(-i,l)} = s^{-1} \, \Big( v^{(-i,l)} - \chi\Big(\frac \nu r \big)e^{-\nu\sqrt{s}} v_\tau^{(-i,l)}\Big), \quad
u_0^{(j,k)} = s^{-1} \, \Big( v^{(j,k)} - \chi\Big(\frac \nu r \big)e^{-\nu\sqrt{s}} v_\tau^{(j,k)}\Big)
\]
of $U^{i,l,\gamma}$ and $U_{M_j}^{(j,k)}$, where $v^{(-i,l)}=-\nabla p_0^{(-i,l)}$ and $v^{(j,k)}=-\nabla p_0^{(j,k)}$. Obviously,
\[
\int_{S_R} \big( p_0^{(-i,l)}\, u_0^{(j,k)} - p_0^{(j,k)}\, u_0^{(-i,l)}\big)\cdot \frac{x}R \, d\sigma = s^{-1} \, (A_R+B_R),
\]
where
\[
A_R = - \int_{S_R} \big( p_0^{(-i,l)}\, \nabla p_0^{(j,k)} - p_0^{(j,k)}\, \nabla p_0^{(-i,l)}\big)\cdot \frac{x}R \, d\sigma
\]
and
\[
B_R = \int_{S_R} \chi\Big(\frac \nu r \big)\, e^{-\nu\sqrt{s}} \, \Big( p_0^{(-i,l)}\, v_\tau^{(j,k)}- p_0^{(j,k)}\, v_\tau^{(-i,l)}\Big)\cdot \frac xR \, d\sigma.
\]
Here,
\begin{eqnarray*}
A_R  & = & - \int_{S_R} \big( p_0^{(-i,l)}\, \partial_r p_0^{(j,k)} - p_0^{(j,k)}\, \partial p_0^{(-i,l)}\big)\, d\sigma \\
& = & -(\mu_i+\mu_j+1)\, R^{\mu_j-\mu_i} \int_\Omega \phi_{i,l}(\omega)\, \phi_{j,k}(\omega)\, d\omega = - (2\mu_j+1)\, \delta_{i,j}\, \delta_{k,l}
\end{eqnarray*}
and
\begin{eqnarray*}
|B_R| & \le & c\, \int_{S_R} \chi\Big(\frac \nu r \big)\, e^{-\nu\, \mbox{\scriptsize Re}\, \sqrt{s}} \, r^{\mu_j-\mu_i-2}\, d\sigma \\
& = & c\, R^{\mu_j-\mu_i} \int_\Omega \chi\big( \nu(\omega)\big)\, e^{-R\, \nu(\omega)\, \mbox{\scriptsize Re}\, \sqrt{s}}\, d\omega
  \le c\, |s|^{-1/2} R^{\mu_j-\mu_i-1},
\end{eqnarray*}
i.~e., $B_R \to 0$ as $R\to \infty$ if $\mu_j<\mu_i+1$. Analogously,
\[
\int_{S_R} \big( P^{i,l,\gamma}\, U_{M_j}^{(j,k)} - P_{M_j}^{(j,k)}\, U^{i,l,\gamma}\big)\cdot \frac{x}R \, d\sigma
  - \int_{S_R} \big( p_0^{(-i,l)}\, u_0^{(j,k)} - p_0^{(j,k)}\, u_0^{(-i,l)}\big)\cdot \frac{x}R \, d\sigma
\]
tends to zero as $R\to \infty$ if $\mu_j<\mu_i+1$. This proves the lemma. \hfill $\Box$ \\

Using the last two lemmas, we can prove the following theorem.

\begin{Th} \label{Ct2}
Suppose that $\mbox{\em Re}\, s \ge 0$, $s\not= 0$, and that $(u,p) \in E_\beta^2(K)\times V_\beta^1(K)$ is a solution of the problem
{\em (\ref{par3})} with the data
\[
f\in E_\beta^0(K)\cap E_\gamma^0(K), \quad g\in X_\beta^1(K)\cap X_\gamma^1(K),\quad\mbox{where } \frac 12 -\lambda_1 < \beta < \frac 12 < \gamma < \frac 32 \, .
\]
If $\gamma-\frac 12$ is not an eigenvalue of the pencil ${\cal N}(\lambda)$, then $(u,p)$ admits the decomposition
\begin{equation} \label{1Ct2}
(u,p) = \eta_s \sum_{j \in J_{\gamma}} \sum_{k=1}^{\sigma_j} c_{j,k}(s)\, \big( u_0^{(-j,k)}, p_0^{(-j,k)}\big) + (v,q),
\end{equation}
where $v\in E_\gamma^2(K)$, $q\in V_\gamma^1(K)$, $J_{\gamma}$ denotes the set of all $j$ such that $0 \le \mu_j < \gamma-\frac 12$ and
\begin{equation} \label{1Ct2b}
c_{j,k}(s) = -\frac{s}{1+2\mu_j} \int_K \big( f(x)\cdot V^{(j,k)}(x,s)+g(x)\, Q^{(j,k)}(x,s) \big)\, dx.
\end{equation}
The remainder $(v,q)$ and the coefficients $c_{j,k}$ satisfy the estimates
\begin{equation} \label{1Ct2a}
\| v\|_{V_\gamma^2(K)} + |s|\, \| v\|_{V_\gamma^0(K)} +  \| q\|_{V_\gamma^1(K)}
\le c\, \Big( \| f\|_{V_\gamma^0(K)} + \| g\|_{V_\gamma^1(K)} + |s| \, \| g\|_{(V_{-\gamma}^1(K))^*}\Big)
\end{equation}
and
\begin{equation} \label{1Ct2c}
\big| c_{j,k}(s)\big|^2 \le c\, |s|^{\gamma-\mu_j-1/2}\, \Big( \| f\|^2_{V_\gamma^0(K)} + \| g\|^2_{V_\gamma^1(K)} + |s|^2 \, \| g\|^2_{(V_{-\gamma}^1(K))^*}\Big)
\end{equation}
with a constant $c$ independent of $s$.
\end{Th}

P r o o f.
If $j\in J_\gamma$ and $\gamma<\frac 32$, then $-1<\gamma-\mu_j-\frac 32 <0$ and, consequently, $M_{j,\gamma}=0$. This means that
\[
U^{j,k,\gamma}(x,s) = u_0^{(-j,k)}(x,s) \quad \mbox{and}\quad P^{j,k,\gamma}(x,s) = p_0^{(-j,k)}(x)=r^{-\mu_j-1}\, \phi_{j,k}(\omega) \ \mbox{ for }j\in J_\gamma.
\]
Furthermore, the number $\beta-\frac 12$ is not an eigenvalue of the pencil ${\cal N}(\lambda)$ since these eigenvalues lie outside the interval $(-1,0)$.
Therefore, the decomposition (\ref{1Ct2}) follows from Lemma \ref{Cl4a}. Let $V^{(j,k)},Q^{(j,k)}$ be the functions (\ref{spec}). We prove the formula (\ref{1Ct2b}). Obviously,
\[
\int_K \big( f\cdot V^{(j,k)}+g\, Q^{(j,k)}\big)\, dx =  a_s\big(\eta_s U_{M_j}^{(j,k)}, \eta_s P_{M_j}^{(j,k)},u,p\big)
  +  a_s\big(v^{(j,k)},q^{(j,k)},u,p\big).
\]
Here $v^{(j,k)}\in E_\delta^2(K)$ and $q^{(j,k)}\in V_\delta^1(K)$ with some $\delta\in \big( \frac 12 -\lambda_1,\frac 12\big)$. From Lemma
\ref{l5b} it follows that $(u,p) \in E_{\beta'}^2(K)\times V_{\beta'}^1(K)$ with arbitrary $\beta'$, $\max( -\delta, \frac 12 -\lambda_1) < \beta'<\frac 12$.
Hence by Lemma \ref{Cl5}, the equality
\[
a_s\big(v^{(j,k)},q^{(j,k)},u,p\big) = a_s\big(u,p,v^{(j,k)},q^{(j,k)}\big) 
\]
holds. Since $a_s\big(u,p,V^{(j,k)},Q^{(j,k)}\big)=0$, it follows that
\begin{eqnarray*}
&& \int_K \big( f\cdot V^{(j,k)}+g\, Q^{(j,k)}\big)\, dx
=  a_s\big(\eta_s U_{M_j}^{(j,k)}, \eta_s P_{M_j}^{(j,k)},u,p\big) - a_s\big(u,p,\eta_s U_{M_j}^{(j,k)},\eta_s P_{M_j}^{(j,k)}\big) \\
&& = A\big(\eta_s U_{M_j}^{(j,k)}, \eta_s P_{M_j}^{(j,k)},u,p\big) .
\end{eqnarray*}
Obviously,
\[
\big|\partial_x^\alpha (\eta_s U_{M_j}^{(j,k)})\big| \le c\, r^{\mu_j-|\alpha|-1}\, |\log r|^{M_j}, \quad
\big|\partial_x^\alpha (\eta_s P_{M_j}^{(j,k)})\big| \le c\, r^{\mu_j-|\alpha|}\, |\log r|^{M_j}
\]
for $|\alpha|\le 2$. Therefore, $\eta_s U_{M_j}^{(j,k)} \in E_{-\gamma}^2(K)$ and $\eta_s P_{M_j}^{(j,k)} \in V_{-\gamma}^1(K)$ if $\mu_j < \gamma-\frac 12$,
and Lemma \ref{Cl5} implies $A\big(\eta_s U_{M_j}^{(j,k)}, \eta_s P_{M_j}^{(j,k)},v,q\big)=0$. Consequently,
\[
\int_K \big( f\cdot V^{(j,k)}+g\, Q^{(j,k)}\big)\, dx = \sum_{i\in J_\gamma} \sum_{l=1}^{\sigma_i}  c_{i,l}(s)\,
    A\big(\eta_s U_{M_j}^{(j,k)}, \eta_s P_{M_j}^{(j,k)}, \eta_s U^{i,l,\gamma},\eta_s P^{i,l,\gamma}\big) .
\]
Using Lemma \ref{Cl6}, we obtain(\ref{1Ct2b}).

We prove the estimates (\ref{1Ct2a}) and (\ref{1Ct2c}). First, let $|s|=1$.
Under the conditions of the theorem, we have $M_1=0$ and $M_j\le 1$ for $0<\mu_j<1$. Since $\mu_j-1$ is not an eigenvalue of the pencil ${\cal N}(\lambda)$ for
$0<\mu_j<1$, we get
\[
\big| \eta_s(r)\, U_{M_j}^{(j,k)}(x,s)\big| \le c\, r^{\mu_j-1}
\]
for $0\le \mu_j < \gamma -\frac 12$ (cf. Remark \ref{Cr1}) and
\begin{equation} \label{2Ct2}
\Big| \int_K \eta_s(r)\, f\cdot U_{M_j}^{(j,k)}(x,s)\, dx\Big|^2 \le c\, \| f\|^2_{V_\gamma^0(K)} \, .
\end{equation}
Analogously, the estimate $\big| \eta_s(r)\, P_{M_j}^{(j,k)}(x,s)\big| \le c\,  r^{\mu_j}$ implies
\begin{equation} \label{3Ct2}
\Big| \int_K \eta_s(r)\, g\cdot P_{M_j}^{(j,k)}(x,s)\, dx\Big|^2 \le c\,  \| g\|^2_{V_{\gamma+1}^0(K)} \, .
\end{equation}
We consider the integral of $f\cdot v^{(j,k)}+g\, q^{(j,k)}$ over $K$. By Lemma \ref{Cl3}, the pair $(v^{(j,k)}, q^{(j,k)})$ is a solution of
the Dirichlet problem for the system
\begin{eqnarray*}
&& (s-\Delta)\, v^{(j,k)} +\nabla q^{(j,k)} =-(s-\Delta)\, \eta_s U_{M_j}^{(j,k)} -\nabla \eta_s P_{M_j}^{(j,k)},  \\
&&  -\nabla\cdot v^{(j,k)} = \nabla\cdot (\eta_s U_{M_j}^{(j,k)})
\end{eqnarray*}
in the space $E_\delta^2(K)\times V_\delta^1(K)$ with some $\delta$ in the interval  $\frac 12 -\lambda_1 < \delta <\frac 12$.
Using Lemma \ref{Cl4}, we get
\begin{eqnarray*}
\big\| (s-\Delta)\, v^{(j,k)} +\nabla q^{(j,k)}\big\|^2_{V_\delta^0(K)}  + \| \nabla\cdot v^{(j,k)}\|^2_{X_\delta^1(K)}
 \le c
\end{eqnarray*}
if $\delta < M_j-\mu_j+\frac 12$. Here, $M_j-\mu_j+\frac 12 > \frac 12 -\lambda_1$. Suppose that
$\frac 12 -\lambda_1 < \delta <\min(M_j-\mu_j+\frac 12,\frac 12)$. Then by Theorem \ref{t3},
\[
\| v^{(j,k)}\|_{E_\delta^2(K)} + \| q^{(j,k)}\|^2_{V_\delta^1(K)}
  \le c.
\]
Since $0< \gamma+\delta< 2$, we have $r^{-2\gamma} \le r^{2\delta-4}$ if $r\le 1$ and
$r^{-2\gamma} \le r^{2\delta}$ if $r\ge 1$. Hence,
\begin{eqnarray*}
\Big| \int_K f\cdot v^{(j,k)}\, dx\Big|^2 & \le & c\, \| f\|^2_{V_\gamma^0(K)}\, \Big( \int_K r^{2\delta-4}\, |v^{(j,k)}|^2\, dx
  + \int_K r^{2\delta}\, |v^{(j,k)}|^2\, dx \\
& \le & c\,   \| f\|^2_{V_\gamma^0(K)}\  \| v^{(j,k)}\|^2_{E_\delta^2(K)}   \le  c'\, \| f\|^2_{V_\gamma^0(K)} \, .
\end{eqnarray*}
Analogously,
\[
\Big| \int_K g\cdot q^{(j,k)}\, dx\Big|^2 \le c\,  \| g\|^2_{X_\gamma^1(K)} \, .
\]
Thus,
\begin{equation} \label{1Ct2d}
\big| c_{j,k}(s)\big|^2 \le c\,  \Big( \| f\|^2_{E_\gamma^0(K)} + \| g\|^2_{X_\gamma^1(K)}\Big) \quad\mbox{for }|s|=1.
\end{equation}
Since the operator $A_\gamma$ is injective for $-\mu_2 -1/2 < \gamma < \lambda_1 +3/2$ (see Lemma \ref{l6}), we get
\[
\| v\|_{E_\gamma^2(K)}  + \| q\|_{V_\gamma^1(K)} \le c\, \Big( \| F\|_{E_\gamma^0(K)}
  + \| G \|_{X_\gamma^1(K)}\Big),
\]
where
\[
F= f - \sum_{j \in J_\gamma} \sum_{k=1}^{\sigma_j} c_{j,k}(s)\, \Big( (s-\Delta)\, (\eta_s U^{j,k,\gamma}) + \nabla (\eta_s P^{j,k,\gamma})\Big)
\]
and
\[
G=g + \sum_{j \in J_{\gamma}} \sum_{k=1}^{\sigma_j} c_{j,k}(s)\, \nabla\cdot (\eta_s U^{j,k,\gamma}) .
\]
Using Lemma \ref{Cl4} and (\ref{1Ct2d}), we get
\[
\| F\|_{E_\gamma^0(K)} + \| G\|_{X_\gamma^1(K)} \le  c\, \Big(  \| f\|_{E_\gamma^0(K)} + \| g\|_{X_\gamma^1(K)}\Big).
\]
This proves (\ref{1Ct2a}) in the case $|s|=1$.

If $s$ is an arbitrary number in the half-plane $\mbox{Re}\, s\ge 0$, $s\not=0$,
we set $x=|s|^{-1/2} y$ and define
\[
\hat{u}(y)=u(x), \ \ \hat{p}(y)=|s|^{-1/2}\, p(x), \ \ \hat{f}(y)=|s|^{-1}\, f(x),\ \ \hat{g}(y)=|s|^{-1/2}\, g(x).
\]
Obviously, $(\hat{u},\hat{p})$ is a solution of the Dirichlet problem for the system
\[
\big( |s|^{-1}s - \Delta\big)\, \hat{u} + \nabla\, \hat{p}=\hat{f}, \ \ -\nabla\cdot \hat{u}=\hat{g}\ \mbox{ in }K.
\]
Consequently,
\[
\big( \hat{u}(y),\hat{p}(y)\big) = \eta(|y|^2) \sum_{j\in J_{\gamma}} \sum_{k=1}^{\sigma_j} c_{j,k}(|s|^{-1}s)\,
  \big( p_0^{(-j,k)}(y), u_0^{(-j,k)}(y,|s|^{-1}s)\big) + \big( \hat{q}(y),\hat{v}(y)\big)
\]
where
\[
p_0^{(-j,k)}(y) =  |s|^{-(\mu_j+1)/2}\, p_0^{(-j,k)}(x), \quad
u_0^{(-j,k)}(y,|s|^{-1}s) = |s|^{-\mu_j/2}\, u_0^{(-j,k)}(x,s)
\]
and
\[
\| \hat{u}\|^2_{E_\gamma^2(K)} + \| \hat{q}\|^2_{V_\gamma^1(K)} + \sum_{j\in J_{\gamma}} \sum_{k=1}^{\sigma_j} \big| c_{j,k}(|s|^{-1}s)\big|^2
\le c\, \Big( \| \hat{f}\|^2_{E_\gamma^0(K)} + \| \hat{g}\|^2_{X_\gamma^1(K)}\Big).
\]
In particular, we get
\[
p(x) = \eta_s(x)  \sum_{j\in J_{\gamma}} \sum_{k=1}^{\sigma_j} |s|^{-\mu_j/2}\, c_{j,k}(|s|^{-1}s)\, p_0^{(-j,k)}(x) + |s|^{1/2}\, \hat{q}(|s|^{1/2}x).
\]
We conclude from this and from (\ref{1Ct2}) that $c_{j,k}(s)= |s|^{-\mu_j/2}\, c_{j,k}(|s|^{-1}s)$ and $q(x)=|s|^{1/2}\, \hat{q}(|s|^{1/2}x)$.
Using the equalities
\[
\| \hat{q}\|^2_{V_\gamma^1(K)} = |s|^{\gamma-1/2}\, \| q\|^2_{V_\gamma^1(K)}, \quad  \| \hat{f}\|^2_{E_\gamma^0(K)} = |s|^{\gamma-1/2}\,
 \| f\|^2_{E_\gamma^0(K)}
\]
and
\[
\| \hat{g}\|^2_{X_\gamma^1(K)} = |s|^{\gamma-1/2} \, \Big( \| g\|^2_{V_\gamma^1(K)} + |s|^2\, \| g\|^2_{(V_{-\gamma}^1(K))^*}\Big),
\]
we obtain
\[
\| q\|^2_{V_\gamma^1(K)} + \sum_{j\in J_{\gamma}} \sum_{k=1}^{\sigma_j} |s|^{\mu_j-\gamma+1/2}\, \big| c_{j,k}(s)\big|^2
  \le c\, \Big( \| f\|^2_{V_\gamma^0(K)} + \| g\|^2_{V_\gamma^1(K)} + |s|^2 \, \| g\|^2_{(V_{-\gamma}^1(K))^*}\Big).
\]
Analogously, the inequality
\[
\|v\|^2_{V_\gamma^2(K)}+ |s|^2\, \|v\|^2_{V_\gamma^0(K)} \le c\, \Big( \| f\|^2_{V_\gamma^0(K)} + \| g\|^2_{V_\gamma^1(K)}
  + |s|^2 \, \| g\|^2_{(V_{-\gamma}^1(K))^*}\Big)
\]
holds. The proof is complete. \hfill $\Box$  \\

The set of all $\mu_j$, $j\in J_\gamma$, contains the simple eigenvalue $\mu_1=0$.
The corresponding singular functions in (\ref{1Ct2}) are
\begin{equation} \label{up0}
p_0^{(-1)}(x)= \frac{|\Omega|^{-1/2}}r \, , \quad u_0^{(-1)}(x,s) = |\Omega|^{-1/2} \Big(\frac{x}{sr^3}
  - \chi\Big( \frac{\nu}r \Big)\, e^{-\nu\sqrt{s}}\,  \Big( \frac{x}{sr^3}- \frac{x\cdot\nabla\nu}{sr^3}\, \nabla\nu\Big)\Big).
\end{equation}
Furthermore, the constant pair $(0,|\Omega|^{-1/2})$ is the solution $(V^{(1)},Q^{(1)})$ in Corollary \ref{Cc1} for the eigenvalue $\mu_1=0$.
Consequently, the coefficient of $(u_0^{(-1)},p_0^{(-1)})$ in Theorem \ref{Ct2} is given by the formula
\begin{equation} \label{coeffup0}
c_1(s) = -s\, |\Omega|^{-1/2} \int_K g\, dx.
\end{equation}
In particular, the term $c_{1}(s) \, ( u_0^{(-1)},p_0^{(-1)})$ does not appear in the asymptotics of $(u,p)$ if
$g$ satisfies the condition (\ref{tildeX}). In this case, the condition on $\gamma$ in the last theorem can be weakened.

\begin{Th} \label{Ct3}
Suppose that $\mbox{\em Re}\, s \ge 0$, $s\not= 0$, and that $(u,p) \in E_\beta^2(K)\times V_\beta^1(K)$ is a solution of the problem
{\em (\ref{par3})} with the data
\[
f\in E_\beta^0(K)\cap E_\gamma^0(K), \quad g\in X_\beta^1(K)\cap X_\gamma^1(K),\quad\mbox{where } \frac 12 -\lambda_1 < \beta < \frac 12 < \gamma
< \min(\lambda_1,\mu_2)+\frac 32 \, .
\]
If $g$ satisfies the condition {\em (\ref{tildeX})} and  $\gamma-\frac 12$ is not an eigenvalue of the pencil ${\cal N}(\lambda)$,
then $(u,p)$ admits the decomposition
\begin{equation} \label{1Ct3}
(u,p) = \eta_s(r) \sum_{j \in J_{\gamma}\backslash \{1\}} \sum_{k=1}^{\sigma_j} c_{j,k}(s)\, \big( u_0^{(-j,k)}, p_0^{(-j,k)}\big) + (v,q),
\end{equation}
where $J_\gamma$ is the same set as in Theorem {\em \ref{Ct2}}, $v\in E_\gamma^2(K)$, $q\in V_\gamma^1(K)$, and the coefficients $c_{j,k}$
are given by the formula {\em (\ref{1Ct2b})}.
The remainder $(v,q)$ and the coefficients $c_{j,k}$ satisfy the estimates {\em (\ref{1Ct2a})} and {\em (\ref{1Ct2c})}
with a constant $c$ independent of $s$.
\end{Th}

P r o o f.
Let $\gamma'$ be an arbitrary real number such that $\gamma'\le \gamma$ and $\frac 12 < \gamma'< \min(\mu_2+\frac 12 , \lambda_1+\frac 32)$.
Then it follows from Lemma \ref{l1} that $u\in E_{\gamma'}^2(K)$ and $p\in V_{\gamma'}^1(K)$.
The interval $\gamma'-\frac 12 < \mu < \gamma-\frac 12$ contains the eigenvalues $\mu_j$ with index $j\in J_\gamma$, $j\not=1$. For these
eigenvalues, the inequalities $-1<\gamma-\mu_j-\frac 32 < 0$ are satisfied. This means that
$M_{j,\gamma}=0$, $U^{j,k,\gamma}(x,s) = u_0^{(-j,k)}(x,s)$ and $P^{j,k,\gamma}(x,s) = p_0^{(-j,k)}(x)$ for $j\in J_\gamma$, $j\not=1$.
Using Lemma \ref{Cl4a}, we obtain (\ref{1Ct3}), where $(v,q) \in E_\gamma^2(K)\times V_\gamma^1(K)$.
Furthermore, as in the proof of Theorem \ref{Ct2}, we obtain
\begin{eqnarray*}
\int_K \Big( f\cdot V^{(j,k)} + g\, Q^{(j,k)} \Big)\, dx = \sum_{i\in J_\gamma\backslash\{ 1\}} \sum_{l=1}^{\sigma_i}  c_{i,l}(s)\,
    A\big(\eta_s U_{M_j}^{(j,k)}, \eta_s P_{M_j}^{(j,k)}, \eta_s U^{i,l,\gamma},\eta_s P^{i,l,\gamma}\big)
\end{eqnarray*}
for $j\in J_\gamma$, $j\not=1$. Obviously, $\mu_j<\mu_i+1$ for $i,j \in J_\gamma\backslash \{1\}$. Therefore, Lemma \ref{Cl6} implies (\ref{1Ct2b})
for $j\in J_\gamma$, $j\not=1$, $k=1,\ldots,\sigma_j$.

Analogously to the proof of Theorem \ref{Ct2}, we obtain (\ref{2Ct2}) and (\ref{3Ct2}) for $|s|=1$.
The functions $v^{(j,k)}$ and $q^{(j,k)}$ belong to the spaces $E_\delta^2(K)$ and  $V_\delta^1(K)$, respectively, where
$\frac 12 -\lambda_1 < \delta< \min(\frac 12 , M_j-\mu_j+\frac 12, 2-\gamma)$. Therefore, analogously to the proof of Theorem \ref{Ct2}, the estimate
\begin{eqnarray*}
\Big| \int_K \big( f\cdot v^{(j,k)} + g\, q^{(j,k)}\big)\, dx\big| \le  c\,   \Big( \| f\|^2_{V_\gamma^0(K)} + \| g\|^2_{X_\gamma^1(K)}  \Big)
\end{eqnarray*}
holds for $|s|=1$. This proves (\ref{1Ct2c}) in the case $|s|=1$. Since the operator $A_\gamma$ is injective by Lemma \ref{l6}, we also obtain the
estimate  (\ref{1Ct2a}) in the case $|s|=1$. In the case $|s|\not=1$, the estimates (\ref{1Ct2a}) and  (\ref{1Ct2c}) can be obtained by means of the
transformation $x=|s|^{-1/2} y$ as in the proof of Theorem \ref{Ct2}. \hfill $\Box$ \\

At the end  of the section, we estimate the functions $V^{(j,k)}$ and $Q^{(j,k)}$
in the formula (\ref{1Ct2b}) for the coefficients $c_{j,k}(s)$, $j\in J_\gamma$. In the case $j=1$ (i.e., $\mu_j=0$), the functions $V^{(j,k)}=0$ and $Q^{(j,k)}=|\Omega|^{-1/2}$ are constant.
Since $\gamma< \min(\lambda_1,\mu_2)+\frac 32$ in Theorems \ref{Ct2} and \ref{Ct3}, it suffices to consider the case $0<\mu_j < \min(\lambda_1,\mu_2)+1$.
Obviously, the number $M_j$ in the definition of $U_{M_j}^{(j,k)}$ and $P_{M_j}^{(j,k)}$ is not greater than 1. More precisely, we have $M_j=0$
if $\mu_j<\lambda_1$ and $M_j=1$ if $\lambda_1\le \mu_j<\lambda_1+1$. Therefore, the functions $U_{M_j}^{(j,k)}$ and $P_{M_j}^{(j,k)}$
contain only logarithmic terms if $\mu_j\ge \lambda_1$ and $\mu_j-1$ is an eigenvalue of the pencil ${\cal N}(\lambda)$ (see Remark \ref{Cr1}).
In the case $0<\mu_j<\mu_2+1$, $\mu_j\not=1$, this is not possible since the intervals $(-1,0)$ and $(0,\mu_2)$ are free of eigenvalue of
the pencil ${\cal N}(\lambda)$.

In the next lemma, we obtain point estimates of the solutions of the problem (\ref{par3}).
Let $N_\beta^{l,\sigma}(K)$ be the weighted H\"older space with the norm
\[
\| u\|_{N_\beta^{l,\sigma}(K)} = \sum_{|\alpha|\le l} \sup_{x\in K} |x|^{\beta-l-\sigma+|\alpha|}\, \big| \partial_x^\alpha u(x)\big|
  + \sum_{|\alpha|=l} \sup_{\substack{x,y\in K \\ 2|x-y|<|x|}} |x|^\beta \, \frac{|\partial_x^\alpha u(x) - \partial_y^\alpha u(y)}{|x-y|^\sigma}\, ,
\]
where $l$ is a nonnegative integer, $\beta$ and $\sigma$ are real numbers, $0<\sigma<1$.

\begin{Le} \label{Cl8}
Let $(u,p) \in E_\delta^2(K) \times V_\delta^1(K)$ be a solution of the problem {\em (\ref{par3})}, where $\frac 12-\lambda_1 < \delta <\frac 12$ and
$f(x)=0$, $g(x)=0$ for $2|s|\, |x|^2<1$. Then the following assertions are true.

{\em 1)} If $2|s|\, |x|^2 <1$, then
\begin{eqnarray*}
\big| \partial_x^\alpha u(x)\big|^2 \le c\, |s|^{\delta+\lambda_1-1/2}\, r^{2(\lambda_1-|\alpha|)}\, (1+|\log |sr^2||)^{2d}\, \| (f,g)\|^2_\delta
  \ \mbox{for }|\alpha|\le 2,\\
\big| \partial_x^\alpha p(x)\big|^2 \le c\, |s|^{\delta+\lambda_1-1/2}\, r^{2(\lambda_1-|\alpha|-1)}\, (1+|\log |sr^2||)^{2d}\, \| (f,g)\|^2_\delta
  \ \mbox{for }|\alpha|\le 1,
\end{eqnarray*}
where the constant $c$ is independent of $f,g,s$. Here
\[
\| (f,g)\|_\delta = \| f\|_{V_\delta^0(K)} + \| g\|_{V_\delta^1(K)} + |s| \, \| g\|_{(V_{-\delta}^1(K))^*}.
\]
Furthermore, $d=1$ only in the case when $\lambda_1=1$ and there exists a generalized eigenvector corresponding to this eigenvalue. Otherwise, $d=0$.

{\em 2)} If  $f\in E_\gamma^0(K) \cap N_{\gamma+\sigma-1/2}^{0,\sigma}(K)$, $g\in X_\gamma^1(K) \cap N_{\gamma+\sigma-1/2}^{1,\sigma}(K)$,
$\gamma\ge \frac 92$ and $2sr^2>1$, then
\begin{eqnarray*}
\big| \partial_x^\alpha u(x)\big|^2 \le c\, |s|^{\gamma+|\alpha|-5/2}\, |x|^{-4}\, ||| (f,g) |||^2_{\gamma,\sigma}\ \mbox{for }|\alpha|\le 2,\\
\big| \partial_x^\alpha p(x)\big|^2 \le c\, |s|^{\gamma-1/2}\, |x|^{-2-2|\alpha|}\, ||| (f,g) |||^2_{\gamma,\sigma}\ \mbox{for }|\alpha|\le 1
\end{eqnarray*}
where
\[
||| (f,g) |||_{\gamma,\sigma} = \| (f,g)\|_\gamma + |s|^{-1}\, \Big( \| f\|_{N_{\gamma+\sigma-1/2}^{0,\sigma}(K)} + \| g\|_{N_{\gamma+\sigma-1/2}^{1,\sigma}(K)}\Big)
\]
and $c$ is independent of $f,g,s$.
\end{Le}

P r o o f. 1) We assume first that $|s|=1$ and that $\chi=\chi(r)$ is a smooth cut-off function equal to one near $r=0$, $\chi(r)=0$ for $2r^2>1$.
Since $f$ and $g$ are zero on the support of $\chi$, we get
\[
-\Delta(\chi u)+ \nabla(\chi p)= F, \ \ \nabla\cdot (\chi u) = u\cdot \nabla\chi \ \mbox{ in K},
\]
where $F= -[\Delta,\chi]\, u + p\, \nabla\chi - s\chi u$. Using the imbedding
\[
V_\delta^2(K) \subset N_{\delta+\sigma-1/2}^{0,\sigma}(K)
\]
(cf. \cite[Lemma 3.6.2]{mr-10}) and regularity results for elliptic systems (see \cite{adn}), we conclude that
\[
\| F\|_{N_{\delta+\sigma-1/2}^{0,\sigma}(K)} +
\| u\cdot \nabla\chi \|_{N_{\delta+\sigma-1/2}^{1,\sigma}(K)} \le c\, \Big( \| u\|_{E_\delta^2(K)} + \| p\|_{V_\beta^1(K)}\Big) \le c'\, \| (f,g)\|_\delta\, .
\]
We may assume that the line $\mbox{Re}\, \lambda = \frac 52 -\delta$ is free of eigenvalues
of the pencil ${\cal L}(\lambda)$. Then it follows from \cite[Theorem 5.2, Corollary 5.1]{mp-78a} that
\[
(\chi u,\chi p) = \sum_{j\in I_\delta} \sum_{k=1}^{\kappa_j} c_{j,k}\, (u_{j,k},p_{j,k}) + (v,q),
\]
where $I_\delta$ is the set of all $j$ such that $0< \mbox{Re}\, \lambda_j < \frac 52-\delta$, $c_{j,k}$ are constants,
$u_{j,k},p_{j,k}$ are of the form
\[
u_{j,k}(x)= r^{\lambda_j}\, \sum_{l=0}^k \frac{(\log r)^l}{l!} \Phi_{j,k-l}(\omega), \ \
p_{j,k}(x) = r^{\lambda_j-1}\, \sum_{l=0}^k \frac{(\log r)^l}{l!} \Psi_{j,k-l}(\omega)
\]
($(\Phi_{j,k},\Psi_{j,k})$ are eigenvectors or generalized eigenvectors of the pencil ${\cal L}(\lambda)$ corresponding to the eigenvalue $\lambda_j$), and
\[
\| v\|_{N_{\delta+\sigma-1/2}^{2,\sigma}(K)}+ \| q\|_{N_{\delta+\sigma-1/2}^{1,\sigma}(K)} + \sum_{j\in I_\delta} |c_{j,k}|
  \le c\, \| (f,g)\|_\delta\, .
\]
If $\lambda_1<1$, then $u_{1,k}$ and $p_{1,k}$ do not contain logarithmic terms. In the case $\lambda_1=1$, the functions $u_{1,k}$,
$p_{1,k}$ contain $\log r$ (with power 1) if there exist generalized eigenvectors corresponding to this eigenvalue
(see \cite[Theorems 5.3.2 and 5.4.1]{kmr-01}). Hence
\[
|c_{j,k}\, u_{j,k}|\le c\, r^{\lambda_1}\, (1+|\log r |^d) \, \| (f,g)\|_\delta\, , \quad
|c_{j,k}\,  p_{j,k}|\le c\, r^{\lambda_1-1} \, (1+|\log r |^d)\, \| (f,g)\|_\delta\,
\]
for $r<1$. Since moreover
\[
\sum_{|\alpha|\le 2} \sup_{x\in K} r^{\delta+|\alpha|-5/2} |\partial_x^\alpha v(x)|
 + \sum_{|\alpha|\le 2} \sup_{x\in K} r^{\delta+|\alpha|-3/2} |\partial_x^\alpha q(x)| \le \, c\, \| (f,g)\|_\delta
\]
and $\delta<\frac 12$, we obtain the desired estimate for $|s|=1$, $2r^2<1$.  If $|s|$ is arbitrary, we set $x=|s|^{-1/2} y$ and define
\[
\hat{u}(y)=u(x), \ \ \hat{p}(y)=|s|^{-1/2}\, p(x), \ \ \hat{f}(y)=|s|^{-1}\, f(x),\ \ \hat{g}(y)=|s|^{-1/2}\, g(x).
\]
Obviously, $(\hat{u},\hat{p})$ is a solution of the Dirichlet problem for the system
\begin{equation} \label{9At1}
\big( |s|^{-1}s - \Delta\big)\, \hat{u} + \nabla\, \hat{p}=\hat{f}, \ \ -\nabla\cdot \hat{u}=\hat{g}\ \mbox{ in }K.
\end{equation}
Since $\hat{f}(y)=0$ and $\hat{g}(y)=0$ for $2|y|^2<1$, we obtain
\begin{eqnarray*}
\big| \partial_y^\alpha \hat{u}(y)\big|^2 & \le & c\, |y|^{2(\lambda_1-|\alpha|)}\, (1+|\log |y|^2|)^{2d}\, \| (\hat{f},\hat{g})\|^2_\delta \\
  & = & c\, |y|^{2(\lambda_1-|\alpha|)}\, (1+|\log |y|^2|)^{2d}\, |s|^{\delta-1/2}\, \| (f,g)\|^2_\delta \ \mbox{ for }|\alpha|\le 2,\ 2|y|^2<1
\end{eqnarray*}
and the analogous estimate for $\partial_y^\alpha \hat{p}(y)$. This proves the assertion 1).

2) We start again with the case $|s|=1$. Since $f(x)$ and $g(x)$ are zero for $2|x|^2<1$,  we have
$\| (f,g)\|_\delta \le c\, \| (f,g)\|_\gamma$ with a certain constant $c$. Thus, by Theorem \ref{t3} and Lemma \ref{Cl4a}, $(u,p)$
admits the decomposition
\[
(u,p) = \eta_s(r) \sum_{j \in J_\gamma} \sum_{k=1}^{\sigma_j} c_{j,k}(s)\, \big( U^{j,k,\gamma}, P^{j,k,\gamma}\big) + (v,q),
\]
where
\[
\| v\|_{E_\gamma^2(K)} + \| q\|_{V_\gamma^1(K)} + \sum_{j\in J_\gamma} |c_{j,k}| \le c\, \| (f,g)\|_\gamma \, .
\]
Since $v(x)=u(x)$ and $q(x)=p(x)$ for $2|x|^2<1$, we get the same estimate for the norm of $(v,q)$ in $V_{\gamma-2}^2(K)\times V_{\gamma-2}^1(K)$.
Obviously,
\[
\big|\partial_x^\alpha \eta_s\, U^{j,k,\gamma}\big| \le c\, r^{-2}\ \mbox{for }|\alpha|\le 2, \quad
\big|\partial_x^\alpha \eta_s\, P^{j,k,\gamma}\big| \le c\, r^{-1-|\alpha|}\ \mbox{for }|\alpha|\le 1
\]
Furthermore, it follows from Lemma \ref{Cl3} that
\[
\| (s-\Delta) v + \nabla q\|_{N_{\gamma+\sigma-1/2}^{0,\sigma}(K)} + \| \nabla\cdot v\|_{N_{\gamma+\sigma-1/2}^{0,\sigma}(K)}
  \le c\, ||| (f,g) |||_{\gamma,\sigma}\, .
\]
Since, moreover,
\[
\| v \|_{N_{\gamma+\sigma-1/2}^{0,\sigma}(K)} \le c\, \| v\|_{V_\gamma^2(K)} \le c\, \, \| (f,g)\|_\gamma \, ,
\]
we conclude from \cite[Theorem 5.1, Corollary 5.1]{mp-78a} that
\[
\| v\|_{N_{\gamma+\sigma-1/2}^{2,\sigma}(K)} + \| q\|_{N_{\gamma+\sigma-1/2}^{1,\sigma}(K)} \le c\, ||| (f,g) |||_{\gamma,\sigma}\, .
\]
In particular,
\[
\big| \partial_x^\alpha v(x)\big| \le c\, r^{-\gamma-|\alpha|+5/2}\, ||| (f,g) |||_{\gamma,\sigma}\ \mbox{for }|\alpha|\le 2
\]
and
\[
\big| \partial_x^\alpha q(x)\big| \le c\, r^{-\gamma-|\alpha|+3/2}\, ||| (f,g) |||_{\gamma,\sigma}\ \mbox{for }|\alpha|\le 1. \quad
\]
This proves the desired estimates for $|s|=1$, $2r^2>1$. In the case $|s|\not=1$, we obtain the estimate analogously to part 1) by means of the coordinate
transformation $x=|s|^{-1/2} y$. \hfill $\Box$ \\

We estimate the functions $V^{(j,k)}$ and $Q^{(j,k)}$ by means of the last lemma.

\begin{Le} \label{Cl7}
Let $\mu_j$ be a positive eigenvalue of the pencil ${\cal N}(\lambda)$, $\mu_j<\min(\lambda_1,\mu_2)+1$, $\mu_j\not=1$.
Then the following estimates are valid for $2|s|\, r^2 <1$:
\begin{eqnarray*}
\big| \partial_x^\alpha V^{(j,k)}(x,s)\big| & \le & c\, |s|^{(|\alpha|-1-\mu_j)/2} \, (|s|^{1/2}r)^{\lambda_1-|\alpha|} \, (1+|\log|sr^2||)^d\
  \mbox{ for }|\alpha|\le 2,\\
\big| \partial_x^\alpha Q^{(j,k)}(x,s)\big| & \le & c\, |s|^{(|\alpha|-\mu_j)/2} \, (|s|^{1/2}r)^{\lambda_1-|\alpha|-1} \, (1+|\log|sr^2||)^d\
  \mbox{ for }|\alpha|\le 1.
\end{eqnarray*}
Here $d$ is the same number as in Lemma {\em \ref{Cl8}}. In the case $2|s|\, r^2 >1$, the estimates
\begin{eqnarray} \label{1Cl7}
\big| \partial_x^\alpha V^{(j,k)}(x,s)\big| & \le & c\, |s|^{(|\alpha|-1-\mu_j)/2} \, (|s|^{1/2}r)^{\mu_j-1}  \mbox{ for }|\alpha|\le 2,\\ \label{2Cl7}
\big| \partial_x^\alpha Q^{(j,k)}(x,s)\big| & \le & c\, |s|^{(|\alpha|-\mu_j)/2} \, (|s|^{1/2}r)^{\mu_j-|\alpha|} \  \mbox{ for }|\alpha|\le 1
\end{eqnarray}
are valid.

The same estimates with an additional factor $(1+|\log|s||)^m$ on the right-hand sides hold if $\mu_j=1$.
\end{Le}

P r o o f. We start with the case $2|s|r^2<1$. In this case, the pair $(V^{(j,k)},Q^{(j,k)})$ coincides with
$(v^{(j,k)},q^{(j,k)})$. Here $(v^{(j,k)},q^{(j,k)})$ is the solution of the problem (\ref{par3}) with the data
\[
f= -(s-\Delta)(\eta_s U_{M_j}^{(j,k)}) + \nabla\, (\eta_s P_{M_j}^{(j,k)}), \quad g= \nabla\cdot (\eta_s U_{M_j}^{(j,k)}),
\]
$v^{(j,k)} \in E_\delta^2(K)$, $q^{(j,k)} \in V_\delta^1(K)$, $\frac 12 -\lambda_1 < \delta < \frac 12$. For $0<\mu_j
< \min(\mu_2,\lambda_1)+1$, $\mu_j\not=1$, the functions $U_{M_j}^{(j,k)}$ and $P_{M_j}^{(j,k)}$ do not contain logarithmic
factors and we get
\[
\| (f,g)\|_\delta \le c\, |s|^{-\delta-\mu_j-1/2}
\]
by means of Lemma \ref{Cl4}. In the case $\mu_j=1$, an additional factor $(1+|\log|s||)^k$ appears on the the right-hand side of the last estimate.
Using Lemma \ref{Cl8}, we obtain the desired estimates for the case $2|s|\, r^2 <1$.

We consider the case $2|s|r^2>1$. Let $\gamma \ge \frac 92$ and let $\gamma-\frac 12$ be not an eigenvalue of the pencil ${\cal N}(\lambda)$.
Furthermore, let $N_j$ be an integer, $N_j>\gamma+\mu_j-\frac 12$. We consider the functions
\[
u^{(j,k)} = V^{(j,k)} - \eta_s U_{N_j}^{(j,k)} = v^{(j,k)} + \eta_s\, \big( U_{M_j}^{(j,k)}-U_{N_j}^{(j,k)}\big)
\]
and $p^{(j,k)} = Q^{(j,k)} - \eta_s P_{N_j}^{(j,k)}$. The pair $(u^{(j,k)},p^{(j,k)})$ is a solution of the problem (\ref{par3}) with the data
\[
f= -(s-\Delta)\, (\eta_s U_{N_j}^{(j,k)}) - \nabla(\eta_s P_{N_j}^{(j,k)}), \quad g = \nabla\cdot (\eta_s U_{N_j}^{(j,k)} )  .
\]
By Lemma \ref{Cl4}, $f\in E_\gamma^0(K)$ and $g\in X_\gamma^1(K)$. Hence, Lemma \ref{Cl8} implies
\begin{eqnarray*}
\big| \partial_x^\alpha u^{(j,k)}\big|^2 \le c\, |s|^{\gamma+|\alpha|-5/2}\, |x|^{-4}\, ||| (f,g) |||^2_{\gamma,\sigma}
\end{eqnarray*}
for $|\alpha|\le 2$, $2|s|r^2>1$. Using Lemma \ref{Cl3}, one can easily show that
\[
||| (f,g) |||^2_{\gamma,\sigma} \le c\, |s|^{-\gamma-\mu_j-1/2}\, \big( 1+|\log|s||\big)^{2m}
\]
(cf. Lemma \ref{Cl4}). This implies
\[
\big| \partial_x^\alpha u^{(j,k)}\big| \le c\, |s|^{(|\alpha|-\mu_j-1)/2}\, \big( 1+|\log|s||\big)^{m}\, (|s| r^2)^{-1}
\]
for $|\alpha|\le 2$, $2|s|r^2>1$. Since
\[
\big| \partial_x^\alpha (\eta_s U_{N_j}^{(j,k)}(x,s))\big|  \le  c\, |s|^{(|\alpha|-1-\mu_j)/2} \, \big( 1+|\log|s||\big)^{m}\, (|s|^{1/2}r)^{\mu_j-1}
\]
for $|\alpha|\le 2$, we obtain
\[
\big| \partial_x^\alpha  V^{(j,k)}(x,s)\big| \le c\, |s|^{(|\alpha|-\mu_j-1)/2} \, \big( 1+|\log|s||\big)^{m}\, (|s|^{1/2}r)^{\mu_j-1}
\]
for $|\alpha|\le 2$, $2|s|r^2>1$. Analogously,
\[
\big| \partial_x^\alpha Q^{(j,k)}\big| \le c\, |s|^{|\alpha|-\mu_j}\, \big( 1+|\log|s||\big)^{m}\, (|s|^{1/2}r)^{\mu_j-|\alpha|}
\]
for $|\alpha|\le 1$, $2|s|r^2>1$. If $\mu_j\not=1$, then these estimates can be improved. As was noted above, the functions
$U_{M_j}^{(j,k)}$ and $P_{M_j}^{(j,k)}$ do not contain logarithmic terms in this case. Thus,
\[
U_{M_j}^{(j,k)}(x,s) = |s|^{-(\mu_j+1)/2}\, U_{M_j}^{(j,k)}(|s|^{1/2}x,|s|^{-1}s), \quad
P_{M_j}^{(j,k)}(x,s) = |s|^{-\mu_j/2}\, P_{M_j}^{(j,k)}(|s|^{1/2}x,|s|^{-1}s)
\]
and, consequently,
\[
V^{(j,k)}(x,s) = |s|^{-(\mu_j+1)/2}\, V^{(j,k)}(|s|^{1/2}x,|s|^{-1}s), \quad
Q^{(j,k)}(x,s) = |s|^{-\mu_j/2}\, Q^{(j,k)}(|s|^{1/2}x,|s|^{-1}s).
\]
Since
\[
\big| \partial_x^\alpha  V^{(j,k)}(x,s)\big| \le c\, r^{\mu_j-1} \ \mbox{for }|\alpha|\le 2, \quad
\big| \partial_x^\alpha  V^{(j,k)}(x,s)\big| \le c\, r^{\mu_j-|\alpha|} \ \mbox{for }|\alpha|\le 1,
\]
if $|s|=1$ and $2r^2>1$, we obtain (\ref{1Cl7}) and (\ref{2Cl7}). The proof is complete. \hfill $\Box$

\section{The time-dependent problem}

In this section, we consider the time-dependent problem (\ref{stokes1}), (\ref{stokes2}). We start with existence and uniqueness theorems
which can be easily deduced from the results in Section 1. Using Theorems \ref{Ct2} and \ref{Ct3}, we describe the behavior of the solutions
at infinity.

\subsection{Solvability results}

Let $Q=K\times {\Bbb R}_+ = K\times (0,\infty)$. By $W_\beta^{2l,l}(Q)$, we denote the weighted Sobolev space of all functions $u=u(x,t)$ on $Q$ with finite norm
\[
\| u\|_{W_\beta^{2l,l}(Q)} = \Big( \int_0^\infty  \sum_{k=0}^l \| \partial_t^k u(\cdot,t) \|^2_{V_\beta^{2l-2k}(K)}\, dt \Big)^{1/2}.
\]
In particular, $W_\beta^{0,0}(Q) = L_2\big({\Bbb R}_+,V_\beta^0(K)\big)$ and
$W_\beta^{2,1}(Q)$ is the set of all $u\in L_2({\Bbb R}_+,V_\beta^2(K))$ such that $\partial_t u \in L_2({\Bbb R}_+,V_\beta^0(K))$.
The space $\stackrel{\circ}{W}\!{}_\beta^{2l,l}(Q)$ is the subspace of all $u\in W_\beta^{2l,l}(Q)$ satisfying
the condition $\partial_t^k u|_{t=0}$ for $x\in K$, $k=0,\ldots,l-1$. Note that
$\partial_t^k u(\cdot,0) \in V_\beta^{2l-2k-1}(K)$ for $u\in W_\beta^{2l,l}(Q)$,
$k=0,\ldots,l-1$ (see \cite[Proposition 3.1]{kozlov-89}).
By \cite[Proposition 3.4]{kozlov-89}, the Laplace transform realizes an isomorphism from $\stackrel{\circ}{W}\!{}_\beta^{2l,l}(Q)$
onto the space $H_\beta^{2l}$ of all holomorphic functions $\tilde{u}(x,s)$ for $\mbox{Re}\, s >0$ with values
in $E_\beta^{2l}(K)$ and finite norm
\[
\| \tilde{u}\|_{H_\beta^{2l}} = \sup_{\gamma>0} \Big( \int_{-\infty}^{+\infty}
  \sum_{k=0}^l |\gamma+i\tau|^{2k}\, \| \tilde{u}(\cdot,\gamma+i\tau)\|^2_{V_\beta^{2l-2k}(K)} \, d\tau  \Big)^{1/2}.
\]
The proof of the analogous result in nonweighted spaces can be found in \cite[Theorem 8.1]{av-64}.

The following theorem proved in \cite[Theorem 3.1]{k/r-16}.

\begin{Th} \label{t6}
Suppose that $f\in L_2\big( {\Bbb R}_+,V_\beta^0(K)\big)$, $g\in L_2\big({\Bbb R}_+,V_\beta^1(K)\big)$,
$\partial_t g\in L_2\big({\Bbb R}_+,(V_{-\beta}^1(K))^*\big)$ and $g(x,0)=0$ for $x\in K$, where
\[
-\lambda_1 +1/2<\beta < \min\big( \mu_2 +1/2\, , \, \lambda_1 +3/2\big), \quad \beta\not= 1/2.
\]
In the case $\beta>1/2$, we assume in addition that
\begin{equation} \label{conditiong}
\int_K g(x,t)\, dx =0 \quad\mbox{for almost all }t.
\end{equation}
Then there exists a uniquely determined solution $(u,p) \in W_\beta^{2,1}(Q) \times L_2\big({\Bbb R}_+, V_\beta^1(K)\big)$
of the problem {\em (\ref{stokes1}), (\ref{stokes2})} satisfying the estimate
\begin{equation} \label{1t6}
\| u\|_{W_\beta^{2,1}(Q)} + \| p\|_{L_2({\Bbb R}_+,V_\beta^1(K))}  \le c\, \Big( \| f\|_{W_\beta^{0,0}(Q)}
  + \| g \|_{L_2({\Bbb R}_+,V_\beta^1(K))} +  \| \partial_t g\|_{L_2({\Bbb R}_+,(V_{-\beta}^1(K))^*)}\Big)
\end{equation}
with a constant $c$ independent of $f$, $g$.
\end{Th}

Analogously, the following theorem can be proved by means of Theorem \ref{Bt1}.

\begin{Th} \label{t7}
Suppose that $f\in L_2\big( {\Bbb R}_+,V_\beta^0(K)\big)$, $g\in L_2\big({\Bbb R}_+,V_\beta^1(K)\big)$,
$\partial_t g\in L_2\big({\Bbb R}_+,(V_{-\beta}^1(K))^*\big)$ and $g(x,0)=0$. Furthermore, we assume that
$\lambda_1=1$, that $\lambda_1$ is a simple eigenvalue of the pencil ${\cal L}(\lambda)$ and that $\beta$
satisfies the inequalities
\[
-\max\big( -\mu_2-1/2,-\mbox{\em Re}\, \lambda_2 +1/2\big) <\beta < \min\big( \mu_2 +1/2\, , \, \mbox{\em Re}\, \lambda_2 +3/2\big), \quad
  \beta\not= \pm 1/2,\  \beta\not=5/2.
\]
In the case $1/2<\beta<5/2$, we assume in addition that $g$ satisfies the condition {\em (\ref{conditiong})}.
Then there exists a uniquely determined solution $(u,p) \in W_\beta^{2,1}(Q) \times L_2\big({\Bbb R}_+, V_\beta^1(K)\big)$
of the problem {\em (\ref{stokes1}), (\ref{stokes2})} satisfying the estimate {\em (\ref{1t6})}.
\end{Th}

Furthermore, the following regularity assertion for the solution can be easily deduced from Lemma \ref{l1} (cf. \cite[Theorem 3.2]{k/r-16}).

\begin{Th} \label{t8}
Suppose that $(u,p) \in W_\beta^{2,1}(Q) \times L_2\big({\Bbb R}_+, V_\beta^1(K)\big)$ is a solution of the problem {\em (\ref{stokes1}), (\ref{stokes2})}
with the data
\begin{equation} \label{E1}
f\in L_2\big( {\Bbb R}_+,V_\beta^0(K) \cap V_\gamma^0(K)\big), \quad
\end{equation}
and
\begin{equation} \label{E2}
\left. \begin{array}{l}
g\in L_2\big({\Bbb R}_+,V_\beta^1(K) \cap V_\gamma^1(K)\big) \\
\partial_t g\in L_2\big({\Bbb R}_+,(V_{-\beta}^1(K))^* \cap (V_{-\gamma}^1(K))^*\big), \end{array} \right\}
\end{equation}
where $\frac 12 -\lambda_1 < \beta,\gamma < \min(\mu_2+\frac 12,\lambda_1+\frac 32)$,
$\beta\not=\frac 12$, $\gamma\not=\frac 12$. In the case $\gamma > \frac 12$, we assume in addition that $g$ satisfies the condition
{\em (\ref{conditiong})}. Then $u\in W_\beta^{2,1}(Q)$ and $p\in L_2\big({\Bbb R}_+, V_\beta^1(K)\big)$.
\end{Th}

If $\lambda_1=1$ and $\lambda_1$ is simple, then this regularity assertion can be improved by means of Lemma \ref{l2}.

\subsection{Asymptotics at infinity}

We consider the solution $(u,p)\in  W_\beta^{2,1}(Q) \times L_2\big({\Bbb R}_+, V_\beta^1(K)\big)$ of the problem (\ref{stokes1}), (\ref{stokes2})
with the data (\ref{E1}), (\ref{E2}). If $\frac 12 - \lambda_1 < \beta < \gamma < \frac 12$, then it follows from Theorem \ref{t8} that
$u\in W_\gamma^{2,1}(Q)$ and $p\in L_2\big({\Bbb R}_+, V_\gamma^1(K)\big)$ .

Now, let $\frac 12 - \lambda_1 < \beta < \frac 12 < \gamma < \frac 32$. We denote the Laplace transforms of $u(x,t)$ and $p(x,t)$
by $\tilde{u}(x,s)$ and $\tilde{p}(x,s)$, respectively. Under the condition of Theorem \ref{Ct2}, we get the decomposition (\ref{1Ct2}) for
$(\tilde{u},\tilde{p})$, i.~e.,
\begin{equation} \label{astilde}
(\tilde{u},\tilde{p}) = \eta_s \sum_{j \in J_{\gamma}} \sum_{k=1}^{\sigma_j} c_{j,k}(s)\, \big( u_0^{(-j,k)}, p_0^{(-j,k)}\big) + (V,Q),
\end{equation}
where $V\in E_\gamma^2(K)$, $Q\in V_\gamma^1(K)$, $J_{\gamma}$ denotes the set of all $j$ such that $0 \le \mu_j < \gamma-\frac 12$ and
\[
c_{j,k}(s) = -\frac{s}{1+2\mu_j} \int_K \big( \tilde{f}(y,s)\cdot V^{(j,k)}(y,s)+\tilde{g}(y,s)\, Q^{(j,k)}(y,s) \big)\, dy.
\]
Let $\psi$ be a $C^\infty$-function on $(-\infty,+\infty)$ with support in the interval $[0,1]$ satisfying the conditions
\[
\int_0^1 \psi(t)\, dt = 1, \quad \int_0^1 t^j\, \psi(t)\, dt = 0 \ \mbox{ for }j=1,\ldots,N,
\]
where $N$ is an integer, $2N> \frac 12 - \gamma-\lambda_1$. By $\tilde{\psi}$, we denote the Laplace transform of $\psi$.
The function $\tilde{\psi}$ is analytic in ${\Bbb C}$ and satisfies the conditions
$\tilde{\psi}(0)=1$, $\tilde{\psi}^{(j)}(0)=0$ for $j=1,\ldots,N$. Since $s^n \tilde{\psi}^{(j)}(s)$ is the Laplace transform of the function
$(-1)^j\, \frac{d^n}{dt^n}\big( t^j\, \psi(t)\big)$, it follows that
\[
\big|\tilde{\psi}^{(j)}(s)\big| \le c_{j,n}\, |s|^{-n}
\]
for every $j\ge 0$ and $n\ge 0$, where $c_{j,n}$ is independent of $s$, $\mbox{Re}\, s\ge 0$. Thus, we can replace the function
$\eta_s(r)=\eta(|s|\, r^2)$ in (\ref{astilde}) by $1-\tilde{\psi}(sr^2)$ and obtain
\begin{eqnarray} \label{D1}
\tilde{u}(x,s) & = & \sum_{j \in J_{\gamma}} \sum_{k=1}^{\sigma_j}  \int_K \big( \tilde{K}_u^{(j,k)}(x,y,s) \, \tilde{f}(y,s)
  + \tilde{H}_u^{(j,k)}(x,y,s)\, \tilde{g}(y,s) \big)\, dy  + V(x,s), \\ \label{D2}
\tilde{p}(x,s) & = & \sum_{j \in J_{\gamma}} \sum_{k=1}^{\sigma_j}  \int_K \big( \tilde{K}_p^{(j,k)}(x,y,s) \cdot \tilde{f}(y,s)
  + \tilde{H}_p^{(j,k)}(x,y,s)\, \tilde{g}(y,s) \big)\, dy   + Q(x,s),
\end{eqnarray}
where
\begin{eqnarray*}
&& \tilde{K}_u^{(j,k)}(x,y,s) = -\frac{s}{1+2\mu_j}\, \big( 1-\tilde{\psi}(sr^2)\big) \, u_0^{(-j,k)}(x,s) \otimes V^{(j,k)}(y,s), \\
&& \tilde{H}_u^{(j,k)}(x,y,s) = -\frac{s}{1+2\mu_j}\, \big( 1-\tilde{\psi}(sr^2)\big) \, Q^{(j,k)}(y,s)\, u_0^{(-j,k)}(x,s)
\end{eqnarray*}
and
\begin{eqnarray*}
&& \tilde{K}_p^{(j,k)}(x,y,s) = -\frac{s}{1+2\mu_j}\, \big( 1-\tilde{\psi}(sr^2)\big) \, p_0^{(-j,k)}(x)\, V^{(j,k)}(y,s), \\
&& \tilde{H}_p^{(j,k)}(x,y,s) = -\frac{s}{1+2\mu_j}\, \big( 1-\tilde{\psi}(sr^2)\big) \, p_0^{(-j,k)}(x)\, Q^{(j,k)}(y,s).
\end{eqnarray*}
By Theorem \ref{Ct3}, the representation (\ref{D1}) holds also if $\gamma<\min(\lambda_1,\mu_2)+\frac 32$ and $g$ satisfies the condition (\ref{conditiong}).
The matrix $\tilde{K}_u^{(j,k)}(x,y,s)$ is the Laplace transform of
\[
K_u^{(j,k)}(x,y,t) = \frac{1}{2\pi i} \, \partial_t^m \int_{-i\infty}^{+i\infty} e^{st}\, s^{-m}\, \tilde{K}_u^{(j,k)}(x,y,s)\, ds
\]
(the integral is absolutely convergent if $m\ge 1$). Analogously, the inverse Laplace transforms $K_p$, $H_u$ and $H_p$ of
$\tilde{K}_p,\tilde{H}_u$ and $\tilde{H}_p$ are defined.
We estimate $K_u^{(j,k)}$, $K_p^{(j,k)}$, $H_u^{(j,k)}$ and $H_p^{(j,k)}$.

\begin{Le} \label{Dl2}
Suppose that $0\le \mu_j < \min(\lambda_1,\mu_2)+1$ and $\mu_j\not=1$. Then
\begin{eqnarray*}
\big| \partial_x^\alpha \partial_y^\beta K_u^{(j,k)}(x,y,t)\big| & \le & c\, t^{-(3+|\alpha|+|\beta|)/2}\, \Big( 1+\frac{|x|}{\sqrt{t}}\Big)^{-2-\mu_j}
  \ \Big( 1 + \frac{|y|}{\sqrt{t}}\Big)^{\mu_j-1}  \\
&& \quad \times \ \Big( \frac{|y|}{|y|+\sqrt{t}}\Big)^{\lambda_1-|\beta|} \ \Big( 1+\log \frac{|y|+\sqrt{t}}{|y|}\Big)^d\
     \ \mbox{ for }|\alpha|,|\beta|\le 2, \\
\big| \partial_x^\alpha H_u^{(j,k)}(x,y,t)\big| & \le & c\, t^{-(4+|\alpha|+|\beta|)/2}\, \Big( 1+\frac{|x|}{\sqrt{t}}\Big)^{-2-\mu_j}
  \ \Big( 1 + \frac{|y|}{\sqrt{t}}\Big)^{\mu_j-|\beta|}   \\ && \quad \times \ \Big( \frac{|y|}{|y|+\sqrt{t}}\Big)^{\lambda_1-|\beta|-1} \
  \Big( 1+\log \frac{|y|+\sqrt{t}}{|y|}\Big)^d  \ \mbox{ for }|\alpha|\le 2, \ |\beta|\le 1.
\end{eqnarray*}
Here, $d$ is the same number as in Lemma {\em \ref{Cl8}}. Furthermore,
\begin{eqnarray*}
\big| \partial_x^\alpha \partial_y^\beta K_p^{(j,k)}(x,y,t)\big| & \le & c\, t^{-(4+|\alpha|+|\beta|)/2}\, \Big( 1+\frac{|x|}{\sqrt{t}}\Big)^{-1-\mu_j-|\alpha|}
  \ \Big( 1 + \frac{|y|}{\sqrt{t}}\Big)^{\mu_j-1}  \\
&& \quad \times \ \Big( \frac{|y|}{|y|+\sqrt{t}}\Big)^{\lambda_1-|\beta|} \ \Big( 1+\log \frac{|y|+\sqrt{t}}{|y|}\Big)^d\
     \ \mbox{ for }|\alpha|\le 1,\ |\beta|\le 2, \\
\big| \partial_x^\alpha H_p^{(j,k)}(x,y,t)\big| & \le & c\, t^{-(5+|\alpha|+|\beta|)/2}\, \Big( 1+\frac{|x|}{\sqrt{t}}\Big)^{-1-\mu_j-|\alpha|}
  \ \Big( 1 + \frac{|y|}{\sqrt{t}}\Big)^{\mu_j-|\beta|}   \\ && \quad \times \ \Big( \frac{|y|}{|y|+\sqrt{t}}\Big)^{\lambda_1-|\beta|-1} \
  \Big( 1+\log \frac{|y|+\sqrt{t}}{|y|}\Big)^d  \ \mbox{ for }|\alpha|, |\beta|\le 1.
\end{eqnarray*}
\end{Le}

P r o o f.
First note that all theorems of this paper are not only valid for $\mbox{Re}\, s \ge 0$ but for all complex $s=\sigma e^{i\theta}$, where $\sigma>0$ and
$|\theta| \le \delta+\frac \pi 2$ with a sufficiently small positive number $\delta$. Therefore, one can replace the path of integration
by the contour $\Gamma_{t,\delta} = \Gamma_{t,\delta}^{(1)} \cup  \Gamma_{t,\delta}^{(2)}$, where
\[
\Gamma_{t,\delta}^{(1)} = \{ s=t^{-1}\, e^{i\theta}: -\delta-\pi/2 < \theta<\delta+\pi/2 \} \quad\mbox{and}\quad
\Gamma_{t,\delta}^{(2)} = \{ s=\sigma e^{\pm i(\delta+\pi/2)}: \ \sigma > t^{-1} \}.
\]
This means, we have
\[
K_u^{(j,k)}(x,y,t) = \frac 1{2\pi i} \int_{\Gamma_{t,\delta}} e^{st}\, \tilde{K}_u^{(j,k)}(x,y,s)\, ds,
\]
Obviously,
\begin{eqnarray*}
&&\big| \partial_x^\alpha (1-\tilde{\psi}(sr^2))\, u_0^{(-j,k)}(x,s)\big| \le c\, |s|^{-1+|\alpha|/2}\, |x|^{-2-\mu_j} \ \mbox{ for }|s|\, r^2\ge 1,\\
&& \big| \partial_x^\alpha (1-\tilde{\psi}(sr^2))\, u_0^{(-j,k)}(x,s)\big| \le c\, |s|^{-1}\, |x|^{-2-\mu_j-|\alpha|}\, |sr^2|^{N+1} \
  \mbox{ for }|s|\, r^2\le 1.
\end{eqnarray*}
If $2N\ge \mu_j+|\alpha|$, this implies
\[
\big| \partial_x^\alpha (1-\tilde{\psi}(sr^2))\, u_0^{(-j,k)}(x,s)\big| \le c\, |s|^{(\mu_j+|\alpha|)/2}\, (1+|s|^{1/2}r)^{-2-\mu_j}.
\]
Using Lemma \ref{Cl7}, we obtain the inequality
\begin{eqnarray*}
\big|  \partial_y^\beta V^{(j,k)}(y,s)\big| & \le & c\, |s|^{(|\beta|-1-\mu_j)/2}\, \Big( \frac{|s|^{1/2}|y|}{1+|s|^{1/2}|y|}\Big)^{\lambda_1-|\beta|} \\
  && \times \ \Big( 1+  \log \frac{1+|s|^{1/2}|y|}{|s|^{1/2}|y|} \Big)^d  \ \big( 1 + |s|^{1/2}|y|\big)^{\mu_j-1} .
\end{eqnarray*}
This directly yields
\begin{eqnarray*}
&& \hspace{-4em}\Big|  \partial_x^\alpha \partial_y^\beta \int_{\Gamma_{t,\delta}^{(1)}} e^{st}\, \tilde{K}_u^{(j,k)}(x,y,s)\, ds\Big|
 \le    c\, t^{-(3+|\alpha|+|\beta|)/2}\, \Big( 1+\frac{|x|}{\sqrt{t}}\Big)^{-2-\mu_j} \\
&& \hspace{6em} \times \ \Big( 1 + \frac{|y|}{\sqrt{t}}\Big)^{\mu_j-1}
  \ \Big( \frac{|y|}{|y|+\sqrt{t}}\Big)^{\lambda_1-|\beta|} \ \Big( 1+\log\frac{|y|+\sqrt{t}}{|y|}\Big)^d.
\end{eqnarray*}
Furthermore,
\begin{eqnarray*}
&& \Big|  \partial_x^\alpha \partial_y^\beta\int_{\Gamma_{t,\delta}^{(2)}} e^{st}\, \tilde{K}_u^{(j,k)}(x,y,s)\, ds\Big|
\le  c\, \Big( 1+ \frac{|x|}{\sqrt{t}}\Big)^{-2-\mu_j} \\
&& \times \int_{1/t}^\infty e^{-\sigma t\sin\delta}\,  \sigma^{(1+|\alpha|+|\beta|)/2} \
   \big( 1 + \sigma^{1/2}|y|\big)^{\mu_j-1}\, \Big( \frac{\sigma^{1/2}|y|}{1+\sigma^{1/2}|y|}\Big)^{\lambda_1-|\beta|} \
   \Big( 1+ \log \frac{1+\sigma^{1/2}|y|}{\sigma^{1/2}|y|}\Big)^d\, d\sigma.
\end{eqnarray*}
In the case $|y|\ge \sqrt{t}$, we have $\sigma^{1/2}|y|\le 1+\sigma^{1/2}|y|\le 2\, \sigma^{1/2}|y|$ for $\sigma\ge 1/t$ and, consequently,
\begin{eqnarray*}
&& \Big|  \partial_x^\alpha \partial_y^\beta \int_{\Gamma_{t,\alpha}^{(2)}} e^{st}\, \tilde{K}_u^{(j,k)}(x,y,s)\, ds\Big|  \\
&& \le   c\, \Big( 1+ \frac{|x|}{\sqrt{t}}\Big)^{-2-\mu_j}\, \int_{1/t}^\infty e^{-\sigma t\sin\delta}\,  \sigma^{(1+|\alpha|+|\beta|)/2} \,
  \big( \sigma^{1/2}|y|\big)^{\mu_j-1}\, d\sigma \\
&& =  c\, t^{-(3+|\alpha|+|\beta|)/2}\, \Big( 1+ \frac{|x|}{\sqrt{t}}\Big)^{-2-\mu_j-|\alpha|} \
  \Big( \frac{|y|}{\sqrt{t}}\Big)^{\mu_j-1}.
\end{eqnarray*}
In the case $|y|\le \sqrt{t}$, we get
\[
\Big|  \partial_x^\alpha \partial_y^\beta \int_{\Gamma_{t,\alpha}^{(2)}} e^{st}\, \tilde{K}_u^{(j,k)}(x,y,s)\, ds\Big|
\le c\, \Big( 1+ \frac{|x|}{\sqrt{t}}\Big)^{-2-\mu_j}\,  (A+B),
\]
where
\[
A  =  \int_{1/t}^{|y|^{-2}} e^{-\sigma t\sin\delta}\, \sigma^{(1+|\alpha|+|\beta|)/2}\, \big(\sigma^{1/2}|y|\big)^{\lambda_1-|\beta|}\,
   \big( 1+|\log (\sigma^{1/2}|y|)|\big)^d\, d\sigma
\]
and
\[
B = \int_{|y|^{-2}}^\infty e^{-\sigma t\sin\delta}\, \sigma^{(1+|\alpha|+|\beta|)/2}\, \big(\sigma^{1/2}|y|\big)^{\mu_j-1}\, d\sigma.
\]
Substituting $\sigma t= \tau$, we get
\begin{eqnarray*}
A  & \le & c\, t^{-(3+|\alpha|+|\beta|)/2} \, \Big( \frac{|y|}{\sqrt{t}}\Big)^{\lambda_1-|\beta|} \int_1^{\infty}  e^{-\tau\sin\delta}\
 \tau^{(1+\lambda_1+|\alpha|)/2}\, \Big( 1+\Big|\log \frac{|y|}{\sqrt{t}}\Big|+\log \tau\Big)^d\,  d\tau \\
& \le &  c\, t^{-(3+|\alpha|+|\beta|)/2} \, \Big( \frac{|y|}{\sqrt{t}}\Big)^{\lambda_1-|\beta|} \ \Big( 1+\Big|\log \frac{|y|}{\sqrt{t}}\Big|\Big)^d.
\end{eqnarray*}
Furthermore,
\begin{eqnarray*}
B & \le & \int_{|y|^{-2}}^\infty e^{-\sigma t\sin\delta}\, \sigma^{(1+|\alpha|+|\beta|)/2}\, \big(\sigma^{1/2}|y|\big)^{\mu_j+1}\, d\sigma \\
  & \le & \int_{1/t}^\infty e^{-\sigma t\sin\alpha}\, \sigma^{(1+|\alpha|+\beta|)/2}\, \big(\sigma^{1/2}|y|\big)^{\mu_j+1}\, d\sigma
   =  c\, t^{-(3+|\alpha|+|\beta|)/2} \, \Big( \frac{|y|}{\sqrt{t}}\Big)^{\mu_j+1}.
\end{eqnarray*}
Hence,
\[
\Big|  \partial_x^\alpha \partial_y^\beta \int_{\Gamma_{t,\delta}^{(2)}} e^{st}\, \tilde{K}_u^{(j,k)}(x,y,s)\, ds\Big|
  \le c\, t^{-(3+|\alpha|+|\beta|)/2}\, \Big( 1+ \frac{|x|}{\sqrt{t}}\Big)^{-2-\mu_j} \ \Big( \frac{|y|}{\sqrt{t}}\Big)^{\lambda_1-|\beta|}\
  \Big( 1+\Big|\log \frac{|y|}{\sqrt{t}}\Big|\Big)^d
\]
if $|y|\le \sqrt{t}$. Thus both in the cases $|y|\ge \sqrt{t}$ and $|y|\le \sqrt{t}$, we obtain the estimate
\begin{eqnarray*}
&& \Big|  \partial_x^\alpha \partial_y^\beta \int_{\Gamma_{t,\delta}^{(2)}} e^{st}\, \tilde{K}_u^{(j,k)}(x,y,s)\, ds\Big|  \\
&& \le   c\, t^{-(3+|\alpha|+|\beta|)/2}\, \Big( 1+\frac{|x|}{\sqrt{t}}\Big)^{-2-\mu_j}\ \Big( 1 + \frac{|y|}{\sqrt{t}}\Big)^{\mu_j-1} \
 \Big( \frac{|y|}{|y|+\sqrt{t}}\Big)^{\lambda_1-|\beta|} \  \Big( 1+\log \frac{|y|+\sqrt{t}}{|y|}\Big)^d .
\end{eqnarray*}
This proves the estimate for the functions $K_u^{(j,k)}$. The other estimates of the lemma can be proved analogously by means of the inequality
\[
\big| \partial_x^\alpha (1-\tilde{\psi}(sr^2))\, p_0^{(-j,k)}(x)\big| \le c\, |s|^{(\mu_j+|\alpha|+1)/2}\, (1+|s|^{1/2}r)^{-1-\mu_j-|\alpha|}.
\]
and the estimate of $Q^{(j,k)}$ in Lemma \ref{Cl7}. \hfill $\Box$ \\

Using Theorem \ref{Ct2} and Theorem \ref{Ct3}, we obtain the following result.

\begin{Th} \label{t9}
Let $(u,p)\in  W_\beta^{2,1}(Q) \times L_2\big({\Bbb R}_+, V_\beta^1(K)\big)$ be a solution of the problem {\em (\ref{stokes1}), (\ref{stokes2})}
with the data {\em (\ref{E1}), (\ref{E2})}.
We assume that $\frac 12 -\lambda_1 < \beta < \frac 12 < \gamma < \frac 32$ and that $\gamma-1/2$ is not an eigenvalue of the pencil
${\cal N}(\lambda)$. Then $u$ and $p$  admit the decomposition
\begin{equation} \label{1t9}
u = \sum_{j \in J_{\gamma}} \sum_{k=1}^{\sigma_j} S^{(j,k)} + v, \quad p= \sum_{j \in J_{\gamma}} \sum_{k=1}^{\sigma_j} T^{(j,k)} + q,
\end{equation}
where $J_{\gamma}$ is the set of all $j$ such that $0 \le \mu_j < \gamma-\frac 12$,
\begin{eqnarray*}
S^{(j,k)}(x,t) & = & \int_0^t \int_K \big( K_u^{(j,k)}(x,y,t-\tau)\, f(y,\tau)  + H_u^{(j,k)}(x,y,t-\tau)\, g(y,\tau)\big)\, dy\, d\tau,\\
T^{(j,k)}(x,t) & = & \int_0^t \int_K \big( K_p^{(j,k)}(x,y,t-\tau)\cdot f(y,\tau)  + H_u^{(j,k)}(x,y,t-\tau)\, g(y,\tau)\big)\, dy\, d\tau.
\end{eqnarray*}
and $(v,q) \in  W_\gamma^{2,1}(Q) \times L_2\big({\Bbb R}_+, V_\gamma^1(K)\big)$.
Here, $K_u^{(j,k)}$ and $H_u^{(j,k)}$ satisfy the estimates of Lemma {\em \ref{Dl2}}. Furthermore, the estimate
\begin{eqnarray*}
&& \| v\|_{W_\gamma^{2,1}(Q)} + \| q\|_{L_2({\Bbb R}_+, V_\gamma^1(K))} \\
&& \le c\, \Big( \| f\|_{L_2\big({\Bbb R}_+,V_\gamma^0(K)\big)}
  + \| g\|_{L_2\big({\Bbb R}_+, V_\gamma^1(K)\big)} + \| \partial_t g\|_{L_2\big({\Bbb R}_+, (V_{-\gamma}^1(K))^*\big)}\Big)
\end{eqnarray*}
is valid with a constant $c$ independent of $f$ and $g$.

The same result holds if $\frac 12 -\lambda_1 < \beta < \frac 12 < \gamma < \min(\lambda_1,\mu_2)+\frac 32$, the number $\gamma-1/2$ is not an eigenvalue of the pencil
${\cal N}(\lambda)$ and $g$ satisfies the condition {\em (\ref{conditiong})}. Then the set $J_\gamma$ in {\em (\ref{1t9})} can be replaced by
$J_\gamma\backslash \{ 1\}$.
\end{Th}

\end{document}